\documentclass[a4paper,12pt]{amsart}
\usepackage{amsmath}
\usepackage{amsthm}
\usepackage{amscd}
\usepackage{amssymb}
\usepackage{amsfonts}
\usepackage{latexsym}
\usepackage[mathscr]{eucal}

\newtheorem{thm}{Theorem}[section]
\newtheorem{prop}[thm]{Proposition}
\newtheorem{lemma}[thm]{Lemma}
\newtheorem{cor}[thm]{Corollary}

\def\graph{\mathop{\rm {graph}}\nolimits}

\def\Im{\mathop{\rm {Im}}\nolimits}

\def\Ker{\mathop{\rm {Ker}}\nolimits}

\begin{document} 

\title{Indecomposable representations of 
quivers on infinite-dimensional Hilbert spaces}
\author{Masatoshi Enomoto}
\address[Masatoshi Enomoto]{College of Business Administration 
and Information Science, 
Koshien University, Takarazuka, Hyogo 665, Japan}      

\author{Yasuo Watatani}
\address[Yasuo Watatani]{Department of Mathematical Sciences, 
Kyushu University, Hakozaki, 
Fukuoka, 812-8581,  Japan}
\maketitle
\begin{abstract}
We study indecomposable representations of 
quivers on separable 
infinite-dimensional Hilbert spaces by bounded operators. 
We consider a complement of Gabriel's theorem for these 
representations. 
Let $\Gamma$ be a finite, connected quiver. 
If its underlying undirected graph 
contains one of extended Dynkin diagrams $\tilde{A_n} \ (n \geq 0)$,
$ \tilde{D_n} \ (n \geq 4)$, 
$\tilde{E_6}$,$\tilde{E_7}$ and $\tilde{E_8}$, then 
there exists an indecomposable representation of $\Gamma$
on separable infinite-dimensional Hilbert spaces.


\medskip\par\noindent
AMS SUBJECT CLASSIFICATION: 46C07, 47A15, 15A21, 16G20, 16G60.

\end{abstract}

\section{Introduction}

We studied  the relative position of
{\it several subspaces } in a separable 
infinite-dimensional Hilbert space in \cite{EW}. In this paper 
we extend it to the relative position of several subspaces 
along quivers. More generally  we study 
representations 
of quivers on infinite-dimensional Hilbert spaces by bounded 
operators.  
We call them Hilbert representations for short. 

Gabriel's theorem says that a connected finite quiver has 
only finitely many indecomposable representations if and only if 
the underlying undirected graph is one of Dynkin diagrams 
$A_n, D_n, E_6, E_7,E_8$ \cite{Ga}. The theory of representations 
of quivers on finite-dimensional vector spaces has 
been developed by  Bernstein-Gelfand-Ponomarev \cite{BGP}, 
Donovan-Freislish \cite{DF}, V. Dlab-Ringel \cite{DR}, 
Gabriel-Roiter \cite{GR}, 
Kac \cite{Ka}, Nazarova \cite{Na} \dots . 

Furthermore locally  scalar representations of quivers 
in the category of Hilbert spaces were introduced by 
Kruglyak and Roiter \cite{KRo}. They associate 
operators and their adjoint operators with arrows 
and classify them up to the unitary 
equivalence.  They proved an 
analog of Gabriel's theorem.  Their study is connected 
with representations of *-algebras generated by 
linearly related orthogonal projections
, see for example, S. Kruglyak, V. Rabanovich and Y. Samoilenko
\cite{KRS}. 

In this paper we study the existence of indecomposable representations of 
quivers on infinite-dimensional Hilbert spaces. We associate 
bounded operators with arrows but we do not associate 
their adjoint operators simultaneously as in \cite{KRo}. 
 
In particular if we consider a certain quiver $\Gamma$ 
 whose  underlying 
undirected graph is the extended 
Dynkin diagram 
$\tilde{D_4}$, 
then indecomposability of  Hilbert representations of $\Gamma$ 
is reduced to  indecomposability of systems of four subspaces 
studied in \cite{EW}. 
We consider 
a complement of Gabriel's theorem for Hilbert 
representations and 
prove one direction: If the underlying undirected 
graph of a finite, connected quiver $\Gamma$ contains
one of extended Dynkin diagrams 
$\tilde{A_n} \ (n \geq 0)$,
$ \tilde{D_n} \ (n \geq 4)$, 
$\tilde{E_6}$,$\tilde{E_7}$ and $\tilde{E_8}$, then 
there exists an indecomposable representation of $\Gamma$
on separable infinite-dimensional Hilbert spaces. The result 
does not depend on the choice of orientation. 
But we cannot prove the converse. 
In fact if the converse were true, then a long standing 
problem in \cite{Ha} on  transitive lattices of subspaces of 
Hilbert spaces would be settled. 

Recall that we study relative position of $n$ subspaces
in a separable  infinite-dimensional Hilbert space in \cite{EW}. 
See Y. P. Moskaleva and Y. S. Samoilenko \cite{MS} on a 
connection with *-algebras generated by projections. 
 Let $H$ be a Hilbert space and $E_1, \dots E_n$ be $n$ subspaces 
in $H$.  Then we say that  ${\mathcal S} = (H;E_1, \dots , E_n)$  
is a system of $n$ subspaces in $H$ or a $n$-subspace system in $H$.
A system ${\mathcal S}$ is called indecomposable if ${\mathcal S}$ 
can not be decomposed into a nontrivial direct sum.  
For any bounded linear operator $A$ on a Hilbert space $K$, we can 
associate a system ${\mathcal S}_A$ of four subspaces in 
$H = K \oplus K$ by 
\[
{\mathcal S}_A = (H;K\oplus 0,0\oplus K,\graph A, \{(x,x) ; x \in K\}).
\]
In particular on a finite dimensional space, Jordan blocks correspond 
to indecomposable systems. Moreover on an infinite dimensional 
Hilbert space, the above system ${\mathcal S}_A$ is indecomposable 
if and only if $A$ is strongly irreducible, which is an 
infinite-dimensional analog of a Jordan block, see
books  by Jiang and Wang \cite{JW},\cite{JW2}. For example, 
a unilateral shift operator is a typical example of strongly 
irreducible operator.  Such a system of four subspaces give an 
indecomposable Hilbert representation of a quiver with underlying
undirected graph $\tilde{D_4}$. We transform these representations 
and make up indecomposable Hilbert representations of other 
quivers in this paper. In finite dimensional case many 
such functors are introduced, see \cite{DF}, for example. 
We follow some of their constructions. But we have not yet proved 
all such functors preserve indecomposability in infinite-dimensional 
Hilbert setting in general.  We have checked the 
indecomposability of the Hilbert representations
in our concrete examples
by our method .

Main theorem of the paper is the following:
Let $\Gamma$ be a finite, connected quiver. 
If its underlying undirected graph 
contains one of extended Dynkin diagrams
 $\tilde{A_n} \ (n \geq 0)$,
$ \tilde{D_n} \ (n \geq 4)$, 
$\tilde{E_6}$,$\tilde{E_7}$ and $\tilde{E_8}$, then 
there exists an indecomposable representation of $\Gamma$
on separable infinite-dimensional Hilbert spaces. There were 
two difficulties which did not appear in finite-dimensional case. 
Firstly we need to find indecomposable, infinite-dimensional
representations of a certain class of $\Gamma$. We constructed them  
by studying 
the relative position of several subspaces along quivers, where 
vertices and arrows are represented by  subspaces and 
 natural inclusion maps. Secondly 
we need to change the orientation of the quiver preserving 
indecomposability. Here comes reflection functors. Being 
different from finite-dimensional case,  
we need to check the co-closedness condition at 
sources to show that indecomposability is preserved under 
reflection functors. 
We introduce a certain nice class, 
called positive-unitary diagonal
Hilbert representations,  such that co-closedness is 
easily checked and preserved under reflection functors 
at any source.

We believe that there exists an analogy between study of 
Hilbert representations of quivers and  
subfactor theory invented by V. Jones \cite{J}. In fact 
Dynkin diagrams also appear in the classification of 
subfactors, see, for example, Goodman, de la Harpe and 
Jones \cite{GHJ}, 
Evans and Kawahigashi \cite{EK}. 
But we have not yet understood the full relations between them.   

There exists a close interplay  between finite-dimensional 
representations of 
quivers and finite-dimensional representations of path algebras 
in purely algebraic sense. Any Hilbert representation  of a quiver 
gives an operator algebra 
representation of the corresponding path algebra. Therefore we expect 
some relation between Hilbert representations of quivers and certain 
operator algebras associated with quivers. There exist some related 
works, see P. Muhly \cite{Mu}, 
D. W. Kribs and S. C. Power \cite{KP} and 
B. Solel \cite{S}.  
But the relation is not so clear for us.

Throughout the paper a projection means an operator 
$e$ with $e^2 = e = e^*$ and 
an idempotent  means an operator $p$ with $p^2 = p$.

In purely algebraic setting, it is known that if a finite-dimensional 
algebra $R$ is not of representation-finite type, then there exist 
indecomposable 
$R$-modules of infinite length as in M. Auslander \cite{Au}. 
Since  we consider bounded operator 
representations on Hilbert spaces, the result in \cite{Au} cannot be applied 
directly.   See a book 
\cite{KR} for infinite length modules. 

The authors are supported by the Grant-in-Aid for Scientific Research of JSPS.

\section{Representations of quivers}

A quiver $\Gamma=(V,E,s,r)$ is a quadruple consisting of 
the set $V$ of vertices, the set $E$ of arrows, 
and two maps $s, r : E \rightarrow V$, which 
associate with each arrow $\alpha \in E$ its 
support $s(\alpha)$ and range  $r(\alpha)$. We 
sometimes denote by $\alpha : x \rightarrow y$ 
an arrow with $x = s(\alpha)$ and $y = r(\alpha)$. 
Thus a quiver is just a directed graph. We 
denote by $|\Gamma|$ 
the underlying undirected graph of a quiver $\Gamma$. 
A quiver $\Gamma$ is said to be connected if $|\Gamma|$ 
is a connected graph. A quiver $\Gamma$ is said to be finite 
if both $V$ and $E$ are finite sets.  
%A path of lengh $m$ 
%is a finite sequence $\alpha = (\alpha_1, \dots,\alpha_m)$ 
%of arrows such that $r(\alpha_k) = s(\alpha_{k+1})$ for 
%$k = 1,\dots, m-1$. Its support is $s(\alpha) = s(\alpha_1)$ 
%and its range is $r(\alpha) = r(\alpha_m)$. A path of 
%length $m \geq 1$ is called a {\it cycle } if its support 
%and range coincide. A cycle of length one is called  a 
%{\it loop}. A quiver is called {\it acyclic} if it contains 
%no cycles. 

\noindent  
{\bf Definition.}
Let $\Gamma=(V,E,s,r)$ be a finite quiver. We say
that $(H,f)$ is a  {\it Hilbert representation} of $\Gamma$ 
if $H=(H_{v})_{v\in V}$  is a family of  Hilbert spaces 
and $f=(f_{\alpha})_{\alpha \in E}$ is a family of
 bounded linear operators $f_{\alpha} : 
H_{s(\alpha)}\rightarrow H_{r(\alpha)}.$

\noindent  
{\bf Definition.} Let $\Gamma=(V,E,s,r)$ be a finite quiver. 
Let $(H,f)$ and $(K,g)$ be Hilbert representations of $\Gamma.$ 
A {\it homomorphism} $T : (H,f) \rightarrow (K,g)$  is a 
family $T =(T_{v})_{v\in V}$ of bounded operators 
$T_v : H_v \rightarrow K_v$ satisfying,  
for any arrow $\alpha \in E$ 
$$
T_{r(\alpha)}f_{\alpha}=g_{\alpha}T_{s(\alpha)}. 
$$
The composition $T \circ S$ of homomorphisms $T$ and $S$ is defined 
by $(T \circ S)_v = T_v \circ S_v$ for $v\in V$.  
Thus we have obtained 
a category $HRep (\Gamma)$ of Hilbert representations of $\Gamma$  

We denote by 
$Hom((H,f),(K,g))$ the set of homomorphisms 
$T : (H,f) \rightarrow (K,g)$. 
We denote by  $End (H,f): =Hom((H,f),(H,f))$  
the set of endomorphisms.  We denote by 
$$
Idem (H,f) : =\{T\in End (H,f) \ | \ T ^2 = T \}
$$
the set of idempotents of $End (H,f)$.  
Let $0 = (0_{v})_{v \in V}$ 
be a family of zero endomorphisms $0_{v}$ and 
$I = (I_{v})_{v\in V}$ 
be a family of identity endomorphisms $I_{v}$. 
The both $0$ and $I$ are in $Idem(H,f)$.

 Let $\Gamma=(V,E,s,r)$ be a finite quiver and 
$(H,f)$, $(W,g)$  be Hilbert representations of $\Gamma.$ 
We say that
$(H,f)$ and $(W,g)$ are {\it isomorphic}, denoted by  
$(H,f)\simeq(W,g)$, 
if there exists an isomorphism $\varphi : (H,f) \rightarrow (W,g)$, 
that is, there exists a family  
$\varphi=(\varphi_{v})_{v\in V}$ of bounded invertible operators
$\varphi_{v}\in B(H_{v},K_{v})$ such that, for any arrow 
$\alpha \in E$, 
$$\varphi_{r(\alpha)}f_{\alpha}=g_{\alpha
}\varphi_{s(\alpha)}.$$

We say that $(H,f)$ is a 
finite-dimensional representation if $H_v$ is 
finite-dimensional for all $v \in V$. 
And $(H,f)$ is an  
infinite-dimensional representation if $H_v$ is 
infinite-dimensional for some $v \in V$. 

\section{Indecomposable representations of quivers}

In this section we  shall introduce a notion of 
indecomposable representation,  
that is, a representation  which cannot be decomposed into a direct 
sum of smaller representations  anymore. 

\bigskip

\noindent  {\bf Definition.}(Direct sum) Let $\Gamma=(V,E,s,r)$ 
be a finite quiver. Let
$(K,g)$ and $(K',g')$ be Hilbert representations
of $\Gamma.$  Define the direct sum 
$(H,f)=(K,g)\oplus(K',g')$  by 
$$
H_{v}=K_{v}\oplus K'_{v}\ (\text{ for } v\in V) 
\ \text{ and } \ 
f_{\alpha}=g_{\alpha }\oplus g'_{\alpha } \ 
(\text{ for } \alpha \in E).
$$

We say that a Hilbert representation $(H,f)$ is zero, denoted by 
 $(H,f) =0 $,  if   
$H_v = 0$ for any $v \in E$. 

\bigskip

\noindent  {\bf Definition.}(Indecomposable representation).  
A Hilbert representation $(H,f)$ of $\Gamma$ 
is called  {\it decomposable} if 
$(H,f)$ is isomorphic to a direct sum of two 
non-zero Hilbert representations.  
A non-zero Hilbert representation $(H,f)$ of $\Gamma$ 
is said to be  {\it indecomposable} if 
it is not decomposable, that is, 
if $(H,f)\cong(K,g) \oplus (K',g')$ 
then $(K,g) \cong 0$ or 
$(K',g') \cong 0$. 

We start with an easy fact. 
Let $H$ be a Hilbert space and $K_1$, $K_2$ be  closed subspaces 
of $H$. Assume that $K_1 \cap K_2 = 0$ and $H = K_1 + K_2$. 
But we do not assume that $K_1$ and $K_2$ are orthogonal. 
Let $T :H \rightarrow H$ be a bounded operator with 
$TK_i \subset K_i$ for $i = 1,2$.  
Define $S_i = T|_{K_i}: K_i \rightarrow K_i$. 
Consider the (orthogonal) direct sum $K_1 \oplus K_2$ 
and the bounded operator $S_1 \oplus S_2$ on $K_1 \oplus K_2$. 
Define a bounded invertible operator 
$\varphi : H \rightarrow K_1 \oplus K_2$ by 
$\varphi(h) = (h_1, h_2)$ for $h = h_1 + h_2$ with 
$h_i \in K_i$, as in the proof of \cite[Lemma 2.1.]{EW}
Then we have 
$T = \varphi^{-1} \circ (S_1 \oplus S_2)\circ  \varphi$.

The following proposition is used frequently 
to show the indecomposability in concrete examples.

\begin{prop}
 Let $(H,f)$be a Hilbert representation of a quiver $\Gamma$.
Then the following conditions are equivalent: 
\begin{enumerate}
\item $(H,f)$ is indecomposable.
\item $ Idem(H,f) = \{0,I\}$. 
\end{enumerate}
\label{prop:indecomposable-idempotent}

\end{prop}
\begin{proof}
$\lnot$(1)$\Longrightarrow$$\lnot$(2):
Assume that $(H,f)$ is not indecomposable.
Then there exist non-zero representations
$(K,g)$ and $(K',g')$ of $\Gamma$, 
such that $(H,f) \cong (K,g) \oplus (K',g').$
For any $x \in V$, define the projection 
$Q_x \in B(K_x\oplus K'_x)$ 
of $K_x\oplus K'_x$ onto $K_x$.  Then $Q := (Q_x)_{x\in V}$ 
is in  $End(K \oplus K', g \oplus g')$,  
because 
\[
Q_{r(\alpha)}(g_{\alpha},g'_{\alpha}) = (g_{\alpha}, 0) 
= (g_{\alpha},g'_{\alpha}) Q_{s(\alpha)}
\] 
for any $\alpha \in E$. 
Since there exists $v, w \in E$ such that 
$K_v \not= 0$ and $K'_w \not= 0$, we have $Q_v \not=0$ 
and $Q_w \not= I$.    
Thus $Q \not=0$ and $Q \not= I$. 
Let $\varphi =(\varphi_x)_{x \in V} : 
(H,f) \rightarrow (K,g) \oplus (K',g')$ 
be an isomorphism. Put $P_x = (\varphi_x)^{-1}Q_x \varphi_x$ 
for $x \in V$ and $P := (P_x)_{x \in x} \in  Idem(H,f)$. Then 
$P\not=0$ and $P\not= I.$

$\lnot$(2)$\Longrightarrow$$\lnot$(1):
Assume that there exists $P \in Idem(H,f)$ with 
$P \not=0$ and $P\not= I$.
Thus 
there exist $v\in V$ \ and $w\in V$ such that 
$P_{v}\neq0_{v}$, $P_{w}\neq I_{w}.$
For any $x \in V$,  define closed subspaces 
\[
K_{x}=P_{x}(H_{x}), \text { and } 
K'_{x} =(I-P_{x})(H_{x}).
\]
Then $K := (K_x)_x \not= 0$, 
$K' := (K'_x)_x \not= 0$ and $H \cong K \oplus K'$.  
For any $\alpha \in E$, let 
$x = s(\alpha)$ and $y = r(\alpha)$. 
Since $f_{\alpha} P_x = P_y f_{\alpha}$, 
we have $f_{\alpha} K_x \subset  K_y$.  
Similarly,  $f_{\alpha}(I - P_x) = (I -P_y) f_{\alpha}$ 
implies that $f_{\alpha} K'_x \subset  K'_y$. 
We can define 
$g_{\alpha} = f_{\alpha}|_{K_x}: K_x \rightarrow K_y$ 
and $g'_{\alpha} = f_{\alpha}|_{K'_x}: K'_x \rightarrow K'_y$. 
Put $g = (g_{\alpha})_{\alpha}$ and $g' = (g'_{\alpha})_{\alpha}$. 
Then $(K,g)$ and $(K',g')$ are representations of $\Gamma$. 
Define $\varphi _x : H_x \rightarrow K_x \oplus K'_x$ by 
$\varphi _x (\xi) = (P_x\xi, (I-P_x)\xi)$ for $\xi \in H_x$. 
Then $\varphi:= (\varphi _x)_{x \in V}
 : (H,f) \rightarrow (K,g) \oplus (K',g')$ is an isomorphism. 
Since $K := (K_x)_x \not= 0$ and 
$K' := (K'_x)_x \not= 0$, 
$(H,f)$ is decomposable.

\end{proof}
\bigskip
\noindent 
{\bf Remark.}(1)
The proof of the above Proposition \ref{prop:indecomposable-idempotent} 
shows that $(H,f)$ is decomposable if and only if there exist non-zero 
families  $K = (K_x)_{x\in V}$ and 
$K = (K'_x)_{x\in V}$ of 
closed subspaces $K_x$ and $K'_x$ of $H_x$ with 
$K_x \cap K'_x = 0$ and $K_x + K'_x = H_x$ such that 
 $f_{\alpha} K_x \subset  K_y$  and 
$f_{\alpha} K'_x \subset  K'_y$ for any arrow  
$\alpha : x \rightarrow y $. 

\noindent
(2)In the statement of the above 
Proposition \ref{prop:indecomposable-idempotent}, we 
cannot replace 
the set  $Idem (H,f)$ of idempotents of endomorphisms 
by the set of projections of endomorphisms. For example, 
let $H_0= \mathbb C ^2$. Fix an angle $\theta$
with $0 < \theta < \pi /2$.  Put $H_1 = \mathbb C(1,0)$ and
$H_2 = \mathbb C(cos\theta, sin\theta)$.  Then the system 
$(H_0;H_1,H_2)$ of two subspaces is isomorphic to 
\[ 
({\mathbb C}^2 ; {\mathbb C}\oplus 0, 0 \oplus {\mathbb C}) 
\cong (\mathbb C; \mathbb C, 0) \oplus (\mathbb C; 0, \mathbb C).
\]
Hence $(H_0;H_1,H_2)$ is decomposable. 
See Example 2 in \cite{EW} and the Remark after it . 
Now consider the following quiver $\Gamma$ : 
\[
\circ_1 \overset{\alpha_1}{\longrightarrow}
 \circ_0 \overset{\alpha_2}\longleftarrow \circ_2
\]
Define a Hilbert 
representation $(H,f)$ of $\Gamma$ by $H = (H_i)_{i= 0,1,2}$ and 
canonical inclusion maps $f_i = f_{\alpha _i} : 
H_i \rightarrow H_0$ 
for $i = 1,2$. Then the Hilbert representation $(H,f)$ is also 
decomposable, see Example 3 below in this paper. 
 But for any $P = (P_i)_{i= 0,1,2} \in End (H,f)$,  
if $P_i \in B(H_i)$ is a projection for $i= 0,1,2$, then  
$P = 0$ or $P = I$. In fact  $P_0 (H_i) \subset H_i$.  
for $i = 1,2$ .  Let $e_1 \in B(H_0)$ and $e_2 \in B(H_0)$ be the
projections of $H_0$ onto $H_1$ and $H_2$.  
Then the $C^*$-algebra 
$C^*(\{e_1,e_2\})$ generated by $e_1$ and $e_2$ is exactly $B(H_0)
\cong M_2(\mathbb C)$. Since $P_0$ commutes with $e_1$ and $e_2$, 
$P_0 = 0$ or  $P_0 = I$. Because 
$P_i = P_0|_{H_i}$, 
 $P_i = 0$ or $P_i = I$ simultaneously. 

\bigskip

\noindent
{\bf Example 1.} Let $\Gamma$ be a loop with one vertex $1$ and 
one arrow $\alpha : 1 \rightarrow 1$, that is, the underlying 
undirected graph is an extended Dynkin 
diagram $\tilde{A_0}$. Let $H_1 =  \ell^2(\mathbb N)$ 
and $f_{\alpha} = S: H_1 \rightarrow H_1$ be a unilateral shift. 
Then  the Hilbert representation $(H,f)$ is infinite-dimensional 
and indecomposable. In fact,  any $T \in Idem(H,f)$ can be identified 
with $T \in  B(\ell^2(\mathbb N))$ with $T^2 = T$ and $TS = ST$.    
Since $T$ commutes with a unilateral shift $S$, the operator 
$T$ is a lower triangular Toeplitz matrix.  Since 
$T$ is an idempotent, $T = 0$ or $T= I$. Thus $(H,f)$ is 
indecomposable.   
Replacing $S$ by $S + \lambda I$ for $\lambda \in {\mathbb C}$, 
we obtain a family of infinite-dimensional, indecomposable  
 Hilbert representations
 $(H^{\lambda}, f^{\lambda})$ of  $\Gamma$. Since 
$(H^{\lambda}, f^{\lambda})$ and $(H^{\mu}, f^{\mu})$ are 
isomorphic if and only if $S + \lambda I$ and $S + \mu I$ 
is similar, we have uncountably  many 
 infinite-dimensional, indecomposable  
 Hilbert representations of  $\Gamma$.

\smallskip
\noindent
{\bf Example 2.} Let $\Gamma=(V,E,s,r)$ be a quiver whose  
underlying undirected graph is an extended Dynkin 
diagram $\tilde{A_n}, \ \ (n \geq 1)$. Then there exist uncountably 
many infinite-dimensional, indecomposable  
 Hilbert representations of  $\Gamma$. For example,  
consider 

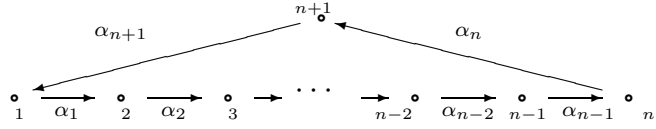
\begin{picture}(300,80)(-30,5)

\put(10,30){\thicklines\circle{2}}
\put(10,25){${}_{{}_{1}}$}
\put(20,30){\vector(1,0){20}}
\put(25,25){${}_{\alpha_{1}}$}

\put(50,30){\thicklines\circle{2}}
\put(50,25){${}_{{}_{2}}$}
\put(60,30){\vector(1,0){20}}
\put(65,25){${}_{\alpha_{2}}$}

\put(90,30){\thicklines\circle{2}}
\put(90,25){${}_{{}_{3}}$}
\put(100,30){\vector(1,0){10}}

\put(115,30){$\cdots$}

\put(140,30){\vector(1,0){10}}
\put(160,30){\thicklines\circle{2}}
\put(145,25){${}_{{}_{n-2}}$}
\put(170,30){\vector(1,0){20}}
\put(170,25){${}_{\alpha_{n-2}}$}

\put(200,30){\thicklines\circle{2}}
\put(195,25){${}_{{}_{n-1}}$}

\put(210,30){\vector(1,0){20}}
\put(215,25){${}_{\alpha_{n-1}}$}

\put(240,30){\thicklines\circle{2}}
\put(245,25){${}_{{}_{n}}$}

\put(125,60){\thicklines\circle{2}}
\put(115,65){${}_{{}_{n+1}}$}
\put(117,58){\vector(-4,-1){100}}
\put(40,55){${}_{\alpha_{n+1}}$}

\put(230,33){\vector(-4,1){100}}
\put(175,55){${}_{\alpha_{n}}$}

\end{picture}

\noindent
Define a Hilbert representation $(H,f)$ of $\Gamma$ by 
$H_1 = H_2 = \dots =H_{n+1} = \ell^2(\mathbb N)$, \ 
$f_{\alpha_2} = f_{\alpha_3} = \dots = f_{\alpha_{n+1}} = I$ 
and $f_{\alpha_{1}} = S$, the unilateral shift. 
Let $P = (P_k)_{k\in V} \in Idem(H,f)$. Then 
\[
P_1 = P_2 = \dots = P_{n+1} \ \text{ and }
SP_1 = P_2S.
\] 
Since $P_1$ is an idempotent and $SP_1 = P_1S$, we have 
$P_1 = 0$ or $P_1 = I$. This implies $P = 0$ or $P = I$. 
Therefore $(H,f)$ is indecomposable. Replacing $S$ by 
$S + \lambda I$ for $\lambda \in {\mathbb C}$, 
we obtain uncountably many 
infinite-dimensional, indecomposable  
 Hilbert representations of  $\Gamma$.  

\smallskip
\noindent
{\bf Example 3.}  
Let $L$ be a Hilbert space and $E_1, \dots E_n$ be $n$ subspaces 
in $L$.  Then we say that  ${\mathcal S} = (L;E_1, \dots , E_n)$  
is a system of $n$ subspaces in $L$. A system ${\mathcal S}$ is 
called indecomposable if ${\mathcal S}$ cannot be decomposed into a 
non-trivial direct sum, see \cite{EW}. 
Consider the following quiver $\Gamma _n=(V,E,s,r)$

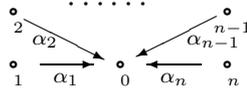
\begin{picture}(100,50)(-100,0)

\put(50,15){\thicklines\circle{2}}
\put(50,10){${}_{{}_{0}}$}

\put(80,15){\vector(-1,0){20}}

\put(65,10){${}_{\alpha_{n}}$}

\put(90,15){\thicklines\circle{2}}
\put(90,10){${}_{{}_{n}}$}

\put(90,35){\thicklines\circle{2}}
\put(85,30){${}_{{}_{n-1}}$}

\put(86,33){\vector(-2,-1){30}}

\put(75,24){${}_{\alpha_{n-1}}$}

\put(10,15){\thicklines\circle{2}}
\put(10,10){${}_{{}_{1}}$}

\put(20,15){\vector(1,0){20}}

\put(25,10){${}_{\alpha_{1}}$}

\put(10,35){\thicklines\circle{2}}
\put(10,30){${}_{{}_{2}}$}

\put(14,32){\vector(2,-1){30}}

\put(17,24){${}_{\alpha_{2}}$}

\put(30,35){$\cdots\cdots$}

\end{picture}

\noindent
Define a Hilbert representation $(H,f)$ of $\Gamma_n$ by 
$H_k := E_k$ \ ($k = 1, \dots, n$), $H_{0} :=L$ and 
$f_k = f_{\alpha_k}  : H_k = E_k \rightarrow  H_{0} =L$ be the 
inclusion map.  Then 
the system ${\mathcal S}$ of $n$ subspaces is 
indecomposable if and only if the corresponding 
Hilbert representation $(H,f)$ of 
$\Gamma$ is indecomposable. 
In fact, assume that  ${\mathcal S}$ is indecomposable. 
Let  $P = (P_k)_{k\in V} \in Idem(H,f)$. Then 
$f_k P_k = P_0f_k$. This implies 
$P_0(H_k) \subset H_k$ for $k = 1, \dots, n$.
Since $P_0$ is idempotent and ${\mathcal S}$ is indecomposable, 
$P_0 = 0$ or $P_0 = I$  by \cite[Lemma 3.2]{EW}. 
Since $f_k P_k = P_0f_k$, $P_k = 0$ or $P_k = I$ simultaneously. 
Thus $P = 0$ or $P = I$, that is, $(H,f)$ is indecomposable. 
Conversely assume that $(H,f)$ is indecomposable.  Let 
$R \in B(L)$ be an idempotent with $R(E_k) \subset E_k$ 
for $k = 1, \dots, n$. Define $P = (P_k)_{k\in V}$  
by $P_0 = R$ and $P_k = P_0|_{H_k}$. Then $P \in Idem(H,f)$. 
Therefore $P = 0$ or $P = I$. Thus $R = O$ or $R = I$. 
Hence ${\mathcal S}$ is indecomposable. 

We can also show that  two systems ${\mathcal S}$ and 
${\mathcal S}'$ of $n$ subspaces are isomorphic 
if and only if the corresponding 
Hilbert representations $(H,f)$  and $(H',f')$  of 
$\Gamma$ are isomorphic. 

Since there exist uncountably many, indecomposable systems 
of fours subspaces in an infinite-dimensional Hilbert space 
as in \cite{EW}, 
there exist uncountably many infinite-dimensional, 
indecomposable Hilbert representations of $\Gamma _4$ 
whose underlying undirected graph is the extended 
Dynkin diagram $\tilde{D_4}$.  

In particular, let $K = \ell^2(\mathbb N)$ and $A \in B(K)$ 
be a strongly irreducible operator studied in \cite{JW}, \cite{JW2} for 
example,  a unilateral shift.  Define \\
$H_0 = K \oplus K$,  $H_1 = K \oplus 0$, 
$H_2 = 0 \oplus K$, \\ 
$H_3 = \{(x,Ax)\in K \oplus K | x \in K\}$, 
$H_4 = \{(x,x)\in K \oplus K | x \in K\}$. \\
Let $f_k = f_{\alpha_k}  : H_k  \rightarrow  H_0 $ be 
the inclusion map for  $k = 1, 2, 3, 4$. 
Put $H^{(A)} = (H_v)_{v \in V}$ and 
$f^{(A)} = (f_{\alpha})_{\alpha \in E}$. 
Then  $(H^{(A)}, f^{(A)})$ is an  
infinite-dimensional, indecomposable  
 Hilbert representation of  $\tilde{D_4}$. 
Moreover let $A$ and $B$ be strongly irreducible 
operators on $\ell^2(\mathbb N)$. Then 
two indecomposable  Hilbert representations 
$(H^{(A)}, f^{(A)})$ and $(H^{(B)}, f^{(B)})$
of  $\tilde{D_4}$ are isomorphic 
if and only if two operators $A$ and $B$ are similar. 

\noindent
{\bf Example 4.}  Consider 
the following quiver $\Gamma =(V,E,s,r)$

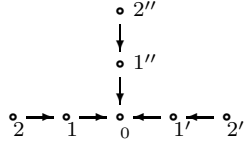
\begin{picture}(130,65)(-80,0)
\put(80,15){\thicklines\circle{2}}
\put(80,10){${}_{{}_{0}}$}

\put(80,30){\vector(0,-1){10}}
\put(80,35){\thicklines\circle{2}}
\put(85,35){${}_{1^{{\prime}{\prime}} }$}

\put(80,50){\vector(0,-1){10}}
\put(80,55){\thicklines\circle{2}}
\put(85,55){${}_{ 2^{{\prime}{\prime}}  }$}

\put(95,15){\vector(-1,0){10}}
\put(100,15){\thicklines\circle{2}}
\put(100,10){${}_{1^{{\prime}} }$}

\put(115,15){\vector(-1,0){10}}
\put(120,15){\thicklines\circle{2}}
\put(120,10){${}_{ 2^{{\prime}}}$}

\put(65,15){\vector(1,0){10}}
\put(60,15){\thicklines\circle{2}}
\put(60,10){${}_{1}$}

\put(45,15){\vector(1,0){10}}
\put(40,15){\thicklines\circle{2}}
\put(40,10){${}_{ 2}  $}

\end{picture}

\noindent
Then underlying undirected graph is an extended Dynkin 
diagram $\tilde{E_6}$. 
Let $K = \ell^2(\mathbb N)$ and $S$ a unilateral shift on $K$. 
We define a Hilbert representation
$(H,f) := ((H_v)_{v\in V},(f_{\alpha})_{\alpha \in E})$ 
of $\Gamma$ as follows: 

Put $H_{0}=K\oplus K\oplus K$, \ $H_{1}=K\oplus0\oplus K$, \ 
$H_{2}=0\oplus0\oplus K$,  \\
$H_{1'}=K\oplus K\oplus0$, \   
$H_{2'}=0\oplus K\oplus0$, \\
$H_{1''}
=\{(x,x,x) + (y,Sy,0) \in K^3\ |\  x, y \in K\}$  and \\
$H_{2''}
=\{(x,x,x) \in K^3  \ | \  x\in K\}$. \\
Then $H_{1''}$ is a closed subspace of $H_0$. In fact, 
let 
\[
(x_n,x_n,x_n ) + (y_n,Sy_n,0) 
= (x_n + y_n,x_n + Sy_n, x_n) \in H_{1''} 
\]
 converges to $(a,b,c) \in H_0$.
Then $x_n \rightarrow c$, 
$y_n \rightarrow a-c$ and $c + S(a-c) = b$. 
Define $x = c$ and $y = a-c$. 
Then $(a,b,c) = (x,x,x) + (y,Sy,0) \in H_{1''}$.  
For any arrow $\alpha \in E$, let  
$f_{\alpha} : H_{s(\alpha)} \rightarrow H_{r(\alpha)}$ be 
the canonical inclusion map. We shall show that 
the Hilbert representation $(H,f)$ is indecomposable. 
Take $T = (T_v)_{v \in V} \in Idem (H,f)$. 
Since $T \in End (H,f)$,  for any 
$v \in \{1,2,1',2',1'', 2'' \}$ and 
any  $x \in H_v$, 
we have $T_0 x = T_vx$. In particular, 
$T_0 H_v \subset H_v$.  
Since 
$H_1\cap H_{1'}=K\oplus0\oplus0$,  
$H_{2'}=0\oplus K\oplus0$ and 
$H_{2}=0\oplus0\oplus K$, 
$T_{0}$ preserves these subspaces. 
Hence $T_{0}$ is a block diagonal operator with
$T_{0}=P\oplus Q\oplus R 
\in B(K \oplus K \oplus K)$. 

Since $T_{0}(H_{2''})\subset H_{2''}$, 
for any $x \in K$, 
\[
T_{0}(x,x,x) = (y,y,y)
\]
for some $y \in K$. Therefore $P=Q=R$ 
and $T_{0}=P\oplus P\oplus P$. 
Moreover $P$ is an idempotent, because so is $T_0$. 
Since $T_0$ preserves 
$H_{1'} \cap H_{1''} =\{(y,Sy,0) \in K^3 \ | \ y\in K \}$,
for any $y \in K$, there exists $z \in K$ such that 
\[
T_{0}\left(
\begin{array}
[c]{c}%
y\\
Sy\\
0
\end{array}
\right)  =\left(
\begin{array}
[c]{c}%
Py\\
PSy\\
0
\end{array}
\right)  =\left(
\begin{array}
[c]{c}%
z\\
Sz\\
0
\end{array}
\right)  .
\]

Therefore $PSy=Sz = SPy$ for any y$\in K$, i.e., $PS=SP$. 
Since $P$ is an idempotent, $P=0$ or $P=I.$
This means that $T_{0}=0$ or $T_{0}=I.$
Because $T_0 x = T_vx$ for any  $x \in H_v$ for 
$v \in \{1,2,1',2',1'', 2'' \}$, 
we have $T_v=0$ or $T_v=I$ simultaneously. 
Thus $T = 0$ or $T = I$, that is, 
$Idem(H,f)=\{0,I\}.$
Therefore $(H,f)$ is indecomposable.

\smallskip
\noindent
{\bf Example 5.} We have a different kind of 
infinite-dimensional, indecomposable  
 Hilbert representation 
$(L,g)= ((L_v)_{v\in V},(g_{\alpha})_{\alpha \in E})$ 
of  the same $\Gamma$ in Example 4 as follows: 
Let $K = \ell^2(\mathbb N)$ and $S$ a unilateral shift on $K$. 
Define $L_{0}=K \oplus K \oplus K,$
$L_{1}=0\oplus K \oplus K,$ \\
$L_{2}=0\oplus\{(y,Sy)\in K^2 \ | \ y\in K\},$
$L_{1^{^{\prime}}}=K \oplus K \oplus0,$ \\
$L_{2^{^{\prime}}}=\{(x,x)\in K^2 \ | \ x\in K\}\oplus0,$
$L_{1^{^{\prime\prime}}}=K\oplus0\oplus K,$ \\
$L_{2^{^{\prime\prime}}}=\{(x,0,x) \in K^3 \ | \ x\in K\}.$
For any arrow $\alpha \in E$, let  
$g_{\alpha} : L_{s(\alpha)} \rightarrow L_{r(\alpha)}$ be 
the canonical inclusion map. We can similarly prove  that 
the Hilbert representation $(L,g)$ is indecomposable. 

We shall show that two Hilbert representations in Example 4 and 5 
are not isomorphic. In fact,  on the contrary, suppose that  
there were an isomorphism 
   $\varphi = (\varphi_v)_{v \in V}: (H,f) \rightarrow (L,g)$. 
Since any arrow is represented by the canonical inclusion, 
 $\varphi_0 : H_0 \rightarrow L_0$ satisfies that 
$\varphi_v = \varphi_0|_{H_v} : H_v \rightarrow L_v$. This 
implies that $\varphi_0(H_v) \subset L_v$ for any $v \in V$. 
Since $\varphi_0(H_{1'}) \subset L_{1'}$ and 
$\varphi_0(H_{1}) \subset L_{1}$, 
$\varphi_0$ has a form such that 
\[
\varphi_{0}=\left(
\begin{array}
[c]{ccc}%
0 & A & 0 \\
B & C & D \\
0 & 0 & E \\
\end{array}
\right).
\]
Since $\varphi_0(H_{2}) \subset L_{2}$, for any $z \in K$ 
there exists $y \in K$ such that $(0,Dz,Ez) = (0,y,Sy)$. 
Hence $Ez = Sy = SDz$, so that $E = SD$. Then 
$\Im \varphi_0 \subset K \oplus K \oplus \Im S \not= L_0$. 
This contradicts the assumption that $\varphi_0 : H_0 \rightarrow L_0$ is 
onto. Therefore Hilbert representations $(H,f)$ and $(L,g)$ 
of $\Gamma$ are not isomorphic.

\section{Reflection functors} 

Reflection functors are crucially used 
in the proof the classification of finite-dimensional, 
indecomposable representations of tame quivers. In fact any 
indecomposable representations of tame quivers can be 
reconstructed by iterating  reflection functors on 
simple indecomposable representations. 
We can not expect such a best situation in infinite-dimensional 
Hilbert representations. But reflection functors are still 
useful to show that some property of representations of 
quivers on infinite-dimensional Hilbert spaces 
does not depend on the choice of orientations 
and does depend on the fact underlying undirected graphs 
are (extended) Dynkin diagrams or not.

Let $\Gamma=(V,E,s,r)$ be a finite quiver. 
A vertex $v \in V$ is called a {\it sink} if $v \not= s(\alpha)$ 
for any $\alpha \in E$. Put 
$E^v = \{ \alpha \in E \ | \ r(\alpha) = v \}$. 
We denote by $\overline{E}$ the set of all 
formally reversed new arrows $\overline{\alpha}$ for 
$\alpha \in E$. 
Thus if $\alpha : x \rightarrow y$ is an arrow, then 
$\overline{\alpha }: x \leftarrow y $. 

\smallskip
\noindent
{\bf Definition.}
Let $\Gamma=(V,E,s,r)$ be a finite quiver.  
For a sink $v \in V$, we construct  a new quiver 
$\sigma_v^+(\Gamma) = (\sigma_v^+(V), \sigma_v^+(E), s,r)$ 
as follows: All the arrows of $\Gamma$ having $v$ as range 
are reversed and all the other arrows remain unchanged.  
More precisely, 
\[
 \sigma_v^+(V) = V \ \ \ 
\sigma_v^+(E) = (E \setminus E^v) \cup \overline{E^v} , 
\]
where $\overline{E^v} = \{ \overline{\alpha } \ | 
\ \alpha \in E^v \}$.  

\smallskip
\noindent
{\bf Definition.} (reflection functor $\Phi_v^+$.) 
Let $\Gamma=(V,E,s,r)$ be a finite quiver.  
For a sink $v \in V$, we define a {\it reflection functor} at  $v$ 
\[
\Phi_v^+ :  HRep (\Gamma) \rightarrow  HRep(\sigma_v^+(\Gamma))
\]  
between the  categories  of Hilbert representations of $\Gamma$  
and $\sigma_v^+(\Gamma)$ as follows: 
For a Hilbert representation $(H,f)$ of $\Gamma$, 
we shall define a Hilbert representation 
$(K,g) = \Phi_v^+(H,f)$ of $\sigma_v^+(\Gamma)$.  
Let
\[
h_v: \oplus_{\alpha \in E^v} H_{s(\alpha)} \rightarrow H_v 
\]
be a bounded linear operator defined by 
\[
h_v((x_{\alpha})_{\alpha \in E^v}) 
= \sum_{\alpha \in E^v} f_{\alpha}(x_{\alpha}).
\]
Define 
\[
K_v := \Ker h_v = \{ (x_{\alpha})_{\alpha \in E^v} \in 
\oplus_{\alpha \in E^v} H_{s(\alpha)}
 \ | \ 
      \sum_{\alpha \in E^v} f_{\alpha}(x_{\alpha}) = 0 \}. 
\]
Consider also the canonical inclusion map 
$i_v : K_v \rightarrow \oplus_{\alpha \in E^v} H_{s(\alpha)}$. 
For $u \in V$ with $u \not=v$, put $K_u = H_u$. 

For $\beta \in E^v$, let 
\[
P_{\beta} : \oplus_{\alpha \in E^v} H_{s(\alpha)} 
\rightarrow H_{s(\beta)}
\]
be the canonical projection. Then define 
\[
g_{\overline{\beta}} :  K_{s(\overline{\beta})} = K_v \rightarrow 
K_{r(\overline{\beta})} = H_{s(\beta)} \ \ 
\text{ by } \  g_{\overline{\beta}} = P_{\beta} \circ i_v
\]
that is, 
$g_{\overline{\beta}}((x_{\alpha})_{\alpha \in E^v}) = x_{\beta}$. 

For $\beta \not\in E^v$, let $g_{\beta} = f_{\beta}$.  

For a homomorphism $T : (H,f) \rightarrow (H',f')$, 
we shall define a homomorphism 
\[
S = (S_u)_{u \in V}= \Phi_v^+(T) : (K,g) 
= \Phi_v^+(H,f) \rightarrow 
(K',g') = \Phi_v^+(H',f')
\] If  $u = v$, a bounded operator 
$S_v : K_v \rightarrow K_v'$ is given by 
\[
S_v((x_{\alpha})_{\alpha \in E^v}) 
= (T_{s(\alpha)}(x_{\alpha}))_{\alpha \in E^v}.
\]

It is easy to see that $S_v$ is well-defined and we have the 
following commutative diagram: 
\[
\begin{CD}
0 @>>> K_v @>i_v>> \oplus_{\alpha \in E^v} H_{s(\alpha)} 
 @>h_v>> H_v \\
@.  @V S_v VV  @V 
(T_{s(\alpha)})_{\alpha \in E^v} VV @V T_v VV \\
0 @>>> {K'}_v  @>{i'}_v>> \oplus_{\alpha \in E^v} {H'}_{s(\alpha)}  
@>h_v'>>{H'}_v
\end{CD}
\]

For other $u \in V$ with $u \not= v$, we put 
\[
S_u = T_u : K_u = H_u \rightarrow K_u' = H_u'. 
\] 

\smallskip

We shall consider a dual of the above construction.  
A vertex $v \in V$ is called a {\it source } 
if $v \not= r(\alpha)$ 
for any $\alpha \in E$. Put 
$E_v = \{ \alpha \in E \ | \ s(\alpha) = v \}$. 

\smallskip
\noindent
{\bf Definition.} 
Let $\Gamma=(V,E,s,r)$ be a finite quiver.  
For a source $v \in V$, we construct  a new quiver 
$\sigma_v^-(\Gamma) = (\sigma_v^-(V), \sigma_v^-(E), s,r)$ 
as follows: All the arrows of $\Gamma$ having $v$ as source 
are reversed and all the other arrows remain unchanged.  
More precisely, 
\[
 \sigma_v^-(V) = V \ \ \ 
\sigma_v^-(E) = (E \setminus E_v) \cup \overline{E_v} , 
\]
where $\overline{E_v} = \{ \overline{\alpha } \ | 
\ \alpha \in E_v \}$.  

In order to define a reflection functor at a source, 
it is convenient to consider the orthogonal complement 
$M^{\perp}$ of a closed subspace $M$ of a Hilbert space $H$ 
instead of the quotient $H/M$. Define an isomorphism 
$f: M^{\perp} \rightarrow H/M$ by $f(y) = [y] = y + M$ 
for $y \in  M^{\perp} \subset H$. Then the inverse 
$f^{-1}:  H/M \rightarrow M^{\perp}$ is given by 
$f^{-1}([x]) = P_M^{\perp}(x)$ for $x \in H$, 
where $P_M^{\perp}$ is the projection of $H$ onto $M^{\perp}$. 
We shall use the following elementary fact frequently:

\begin{lemma}
Let $K$ and $L$ be Hilbert spaces, $M \subset K$ and 
$N \subset L$ be closed subspaces. Let $A: K \rightarrow L$ 
be a bounded operator.  Assume that $A(M) \subset N$.  
Let $\tilde{A} : K/M \rightarrow L/N$ be the induced 
map such that $\tilde{A}([x]) = [Ax]$ for $x \in K$. 
Identifying $K/M$ and $L/N$ with $M^{\perp}$ and $N^{\perp}$,
$\tilde{A}$ is identified with  the bounded operator 
$S: M^{\perp} \rightarrow N^{\perp}$ such that 
$S(x) = P_N^{\perp}(Ax)$.  Then $S = (A^*|_{N^{\perp}})^*$. 
\label{lemma:replacing-quotient}
\end{lemma}
\begin{proof}
Consider $A^* : L \rightarrow K$. Since $A(M) \subset N$, 
we have $A^*(N^{\perp}) \subset M^{\perp}$. Hence 
the restriction $A^*|_{N^{\perp}} : 
N^{\perp} \rightarrow M^{\perp}$ has the adjoint 
\[
(A^*|_{N^{\perp}})^* : M^{\perp} \rightarrow N^{\perp}.
\]
For any $m \in M^{\perp}$ and $n \in N^{\perp}$
\[
((A^*|_{N^{\perp}})^*m|n) = (m|A^*|_{N^{\perp}}n) 
= (m|A^*n) = (Am|n) = (P_N^{\perp}(Am)|n). 
\]

\end{proof}
\smallskip
\noindent
{\bf Definition.} (reflection functor $\Phi_v^-$.) 
Let $\Gamma=(V,E,s,r)$ be a finite quiver.  
For a source $v \in V$, we define a {\it reflection functor} at $v$ 
\[
\Phi_v^- :  HRep (\Gamma) \rightarrow  HRep(\sigma_v^-(\Gamma))
\]  
between the  categories  of Hilbert representations of $\Gamma$  
and $\sigma_v^-(\Gamma)$ as follows: 
For a Hilbert representation $(H,f)$ of $\Gamma$, 
we shall define a Hilbert representation 
$(K,g) = \Phi_v^-(H,f)$ of $\sigma_v^-(\Gamma)$.  
Let
\[
\hat{h}_v: H_v  \rightarrow \oplus_{\alpha \in E_v} H_{r(\alpha)}   
\]
be a bounded linear operator defined by 
\[
\hat{h}_v(x) 
= (f_{\alpha}(x))_{\alpha \in E_v} \ \text{ for } x \in H_v.
\]
Define 
\[
K_v := (\Im \hat{h}_v)^{\perp} = \Ker \hat{h}_v^* \subset 
\oplus_{\alpha \in E_v} H_{r(\alpha)},  
\]
where $\hat{h}_v^* : \oplus_{\alpha \in E_v} H_{r(\alpha)} 
\rightarrow H_v$ is given 
$\hat{h}_v^*((x_{\alpha})_{\alpha \in E_v}) 
= \sum f_{\alpha}^*(x_{\alpha})$. 
For $u \in E$ with $u \not=v$, put $K_u = H_u$. 

Let $Q_v : \oplus_{\alpha \in E_v} H_{r(\alpha)} \rightarrow K_v$ 
be the canonical projection.  
For $\beta \in E_v$, let 
\[
j_{\beta} : H_{r(\beta)} 
\rightarrow \oplus_{\alpha \in E_v} H_{r(\alpha)}
\]
be the canonical inclusion. 
Define 
\[
g_{\overline{\beta}} :  K_{s(\overline{\beta})} = 
H_{r(\beta)} \rightarrow K_{r(\overline{\beta})} = K_v \ \ 
\text{ by }  g_{\overline{\beta}} = Q_v \circ j_{\beta} .
\]

For $\beta \not\in E_v$, let $g_{\beta} = f_{\beta}$.  

For a homomorphism $T : (H,f) \rightarrow (H',f')$, 
we shall define a homomorphism \[
S = (S_u)_{u \in V}= \Phi_v^-(T) : (K,g) 
= \Phi_v^-(H,f) \rightarrow 
(K',g') = \Phi_v^-(H',f'),  
\]
recalling the above Lemma \ref{lemma:replacing-quotient}. 
For $u = v$, a bounded operator 
$S_v : K_v \rightarrow K_v'$ is given by 
\[
S_v((x_{\alpha})_{\alpha \in E_v}) 
= Q_v'((T_{r(\alpha)}(x_{\alpha}))_{\alpha \in E_v}),
\]
where $Q_v' : \oplus_{\alpha \in E_v} H'_{r(\alpha)} 
\rightarrow K'_v$ 
be the canonical projection. 

We have the 
following commutative diagram: 

\[
\begin{CD}
H_v @>\hat{h}_v>> \oplus_{\alpha \in E_v} H_{r(\alpha)} 
@>Q_v>> K_v @>>> 0 \\
@V T_v VV  @V \oplus_{\alpha \in E_v} T_{r(\alpha)} VV  @V S_v VV @.\\
{H'}_v @>{{\hat{h}}'}_v>> \oplus_{\alpha \in E_v} {H'}_{r(\alpha)} 
 @>{Q'}_v>> {K'}_v
@>>>0
\end{CD}
\]
For other $u \in V$ with $u \not= v$, we put 
\[
S_u = T_u : K_u = H_u \rightarrow K_u' = H_u'. 
\] 

We shall explain a relation between two (covariant) functors 
$\Phi_v^+$ and $\Phi_v^-$. We need to introduce another 
(contravariant) functor $\Phi^*$ in the first place.

Let $\Gamma=(V,E,s,r)$ be a finite quiver.  We define 
the opposite quiver 
$\overline{\Gamma}=(\overline{V},\overline{E},s,r)$ 
by reversing all the arrows, that is, 
\[
 \overline{V} = V \ \ \text{ and  } \ \ 
\overline{E} = \{ \overline{\alpha} \ | \ \alpha \in E \}.
\]
\smallskip
\noindent

{\bf Definition.} Let $\Gamma=(V,E,s,r)$ be a finite quiver 
and $\overline{\Gamma}=(\overline{V},\overline{E},s,r)$ 
its opposite quiver. 
We introduce a contravariant functor   
\[
\Phi^*  :  HRep (\Gamma) \rightarrow  HRep(\overline{\Gamma})
\]  
between the  categories  of Hilbert representations of $\Gamma$  
and $\overline{\Gamma}$ as follows: 
For a Hilbert representation $(H,f)$ of $\Gamma$, 
we shall define a Hilbert representation 
$(K,g) = \Phi^*(H,f)$ of $\overline{\Gamma}$ by 
\[
K_u = H_u  \text{ for } u \in V \ \ \ \text{ and } 
g_{\overline{\alpha}} = f_{\alpha}^* \text{ for } 
\alpha  \in E.
\]
For a homomorphism $T : (H,f) \rightarrow (H',f')$, 
we shall define a homomorphism 
\[
S = (S_u)_{u \in V} = \Phi^*(T) : (K',g') 
= \Phi^*(H',f') \rightarrow 
(K,g) = \Phi^*(H,f),
\]  by bounded operators 
$S_u : K_u'= H_u' \rightarrow K_u = H_u$ given by 
$S_u = T_u^*$. 

\begin{prop}
Let $\Gamma=(V,E,s,r)$ be a finite quiver. 
If $v \in V$ is a source of $\Gamma$, then 
$v$ is a sink of $\overline{\Gamma}$, 
$\sigma_v^{-}(\Gamma)
=\overline{\sigma_v^{+}(\overline{\Gamma})}$ 
and we have the following: 
\begin{itemize}
\item [(1)] 
For a  Hilbert representation $(H,f)$ of $\Gamma$, 
\[
\Phi_v^{-}(H,f)=\Phi^* (\Phi_v^{+}(\Phi^*(H,f))).
\]
\item[(2)]
For a homomorphism $T : (H,f) \rightarrow (H',f')$, 
\[
\Phi_v^{-}(T)=\Phi^* (\Phi_v^{+}(\Phi^*(T))).
\]
\end{itemize}
\label{prop:+=*-*}
\end{prop}
\begin{proof}(1):
It is enough to consider around a source $v$. 
For each   $\alpha \in E_v$ with 
$\alpha : v \rightarrow u = r(\alpha)$, 
a bounded operator $f_{\alpha}: H_v \rightarrow H_u$ is 
assigned in $(H,f)$. Taking  $\Phi^*$, we have 
$\Phi^*(H_u) = H_u$  and 
$\Phi^*(f_{\alpha}) = f_{\alpha}^* : H_u \rightarrow H_v$ 
in $\Phi^*(H,f)$. 
Let 
\[
h_v: \oplus_{\alpha \in E_v} H_{r(\alpha)} \rightarrow H_v 
\]
be a bounded operator given  by 
\[
h_v((x_{\alpha})_{\alpha \in E_v}) 
= \sum_{\alpha \in E_v} f_{\alpha}^*(x_{\alpha}).
\]
Define 
\[
W_v := \{ (x_{\alpha})_{\alpha \in E_v} \in 
\oplus_{\alpha \in E_v} H_{r(\alpha)}
 \ | \ 
      \sum_{\alpha \in E_v} f_{\alpha}^*(x_{\alpha}) = 0 \}.
\]
Then $\Phi^{+}_v(\Phi^*(H_v)) = W_v$  and 
$\Phi^{+}_v(\Phi^*(H_u)) = H_u$  
in $\Phi^{+}(\Phi^*(H,f))$. 
Consider the canonical inclusion map 
$i_v : W_v \rightarrow \oplus_{\alpha \in E_v} H_{r(\alpha)}$. 
For $\beta \in E_v$, let 
\[
P_{\beta} : \oplus_{\alpha \in E_v} H_{r(\alpha)} 
\rightarrow H_{r(\beta)}
\]
be the canonical projection.
Then  $\Phi^{+}_v(\Phi^*(f_{\beta})) = P_{\beta} \circ i_v$. 
Finally take $\Phi^*$ again. 
Since $h_v^*: H_v \rightarrow 
\oplus_{\alpha \in E_v} H_{r(\alpha)}$ 
is given by 
\[
(h_v^*)(y) = (f_{\alpha}(y))_{\alpha \in E_v} = \hat{h}_v(y), 
\ \ \text{ for } \ y \in H_v.
\] 
we have 
\[
\Phi^*(\Phi^{+}_v(\Phi^*(H_v))) = W_v
 = \Ker h_v = (\Im h_v^*)^{\perp} 
= (\Im \hat{h}_v)^{\perp}= \Phi^{-}_v(H_v).
\]
Moreover
$i_v^* = Q_v: 
\oplus_{\alpha \in E_v} H_{r(\alpha)} \rightarrow W_v $ 
is the canonical projection. 
For $\beta \in E_v$, we have 
\[
P_{\beta}^* = j_{\beta} : 
H_{r(\beta)} \rightarrow \oplus_{\alpha \in E_v} H_{r(\alpha)} .
\]
Therefore 
\[
\Phi^*(\Phi^{+}_v(\Phi^*(f_{\beta})))
= (P_{\beta} \circ i_v)^* = i_v^* \circ P_{\beta}^* 
= Q_v \circ j_{\beta} 
= \Phi^{-}_v(f_{\beta}). 
\]

\noindent
(2): If $u \not= v$, then 
\[
(\Phi^* (\Phi_v^{+}(\Phi^*(T))))_u  = T_u^{**}= T_u = (\Phi_v^{-}(T))_u .
\]
If $u = v$, then, apply 
Lemma \ref{lemma:replacing-quotient} by putting that 
$K = \oplus_{\alpha \in E_v} H_{r(\alpha)}$, 
$L = \oplus_{\alpha \in E_v} H_{r(\alpha)}'$, 
$M$ is the closure of $ \{ (f_{\alpha}(x))_{\alpha \in E_v} \in K \ | \ x \in H_v \}$ in $K$, 
$N$ is the closure of $\{ (f_{\alpha}'(x))_{\alpha \in E_v} \in L \ | \ x \in H_v' \}$ in $L$
and $A: K \rightarrow L$ with 
$A((y_{\alpha})_{\alpha \in E_v}) 
= (T_{r(\alpha)}y_{\alpha})_{\alpha \in E_v}$. 
Then 
$(\Phi^* (\Phi_v^{+}(\Phi^*(T))))_v  =  (\Phi_v^{-}(T))_v.$

\end{proof}

\begin{prop}
Let $\Gamma=(V,E,s,r)$ be a finite quiver. 
If $v \in V$ is a sink of $\Gamma$, then 
$v$ is a source of $\overline{\Gamma}$, 
$\sigma_v^{+}(\Gamma)
=\overline{\sigma_v^{-}(\overline{\Gamma})}$ 
and we have the following: 
\begin{itemize}
\item [(1)] 
For a  Hilbert representation $(H,f)$ of $\Gamma$, 
\[
\Phi_v^{+}(H,f)=\Phi^* (\Phi_v^{-}(\Phi^*(H,f))). 
\]
\item[(2)]
For a homomorphism $T : (H,f) \rightarrow (H',f')$, 
\[
\Phi_v^{+}(T)=\Phi^* (\Phi_v^{-}(\Phi^*(T))).
\]
\end{itemize}
\end{prop}

\begin{proof}
It follows immediately from Proposition \ref{prop:+=*-*}
and the fact that $(\Phi^*)^2 = Id$.  

\end{proof}

\section{Duality theorem}

We shall show a certain duality between reflection functors.
Bernstein-Gelfand-Ponomarev \cite{BGP} 
introduced reflection functors and Coxeter functors 
and clarify a relation with the Coxeter-Weyl group and 
Dynkin diagrams in the case of finite-dimensional 
representations of quivers. 
In the case of infinite-dimensional Hilbert representations, 
duality theorem 
between reflection functors does not hold as in the 
purely algebraic setting. We need to modify and assume a certain 
closedness condition at a sink or a source. 

\noindent  {\bf Definition.} Let $\Gamma=(V,E,s,r)$ be a finite quiver and 
 $v \in V$ a sink.  Recall that 
$E^v = \{\alpha \ | \ r(\alpha) = v\}$. 
We say that a Hilbert representation $(H,f)$  of $\Gamma$ 
is {\it closed} at $v$ if 
$\sum_{\alpha \in E^v} \Im f_{\alpha} \subset H_v$ 
is a closed subspace.  We say that 
$(H,f)$ is {\it full} at $v$ if 
$\sum_{\alpha \in E^v} \Im f_{\alpha} = H_v$. 

\smallskip

\noindent {\bf Remark.} Recall that a bounded operator 
$h_v: \oplus_{\alpha \in E^v} H_{s(\alpha)} \rightarrow H_v$ 
is given by 
$h_v((x_{\alpha})_{\alpha \in E^v}) 
= \sum_{\alpha \in E^v} f_{\alpha}(x_{\alpha}).$
Then  a Hilbert representation $(H,f)$  of $\Gamma$ 
is {\it closed} at $v$ if and only if  $\Im h_v$ is closed. 
A Hilbert representation $(H,f)$ is {\it full} at $v$ 
if and only if $h_v$ is onto. 

\noindent {\bf Definition.} Let $\Gamma=(V,E,s,r)$ be a finite quiver and 
 $v \in V$ a source.  Recall that $E_v = \{\alpha | s(\alpha) = v\}$.  
We say that a Hilbert representation $(H,f)$  of $\Gamma$ 
is {\it co-closed} at $v$ if 
$\sum_{\alpha \in E_v} \Im f_{\alpha}^* \subset H_v$ 
is a closed subspace.  We say that 
$(H,f)$ is {\it co-full} at $v$ if 
$\sum_{\alpha \in E_v} \Im f_{\alpha}^* = H_v$. 

\smallskip

\noindent 
{\bf Remark.} Recall that a bounded operator 
$
\hat{h}_v: H_v  \rightarrow \oplus_{\alpha \in E_v} H_{r(\alpha)}   
$
is given by 
$
\hat{h}_v(x) 
= (f_{\alpha}(x))_{\alpha \in E_v} \ \text{ for } x \in H_v.
$
Then  a Hilbert representation $(H,f)$  of $\Gamma$ 
is co-closed at $v$ if and only if  $\Im \hat{h}_v^*$ is closed. 
A Hilbert representation $(H,f)$ is  co-full at $v$ 
if and only if $\hat{h}_v^*$ is onto if and only if 
$\Im \hat{h}_v$ is closed and 
$\cap_{\alpha \in E_v} \Ker f_{\alpha} = 0$. 
In fact the latter condition is equivalent to 
$(\Im \hat{h}_v^*)^{\perp}= \Ker \hat{h}_v = 0$. 
We also see that $(H,f)$  
is co-closed at $v$ if and only if $\Phi_v^*(H,f)$ is 
closed at $v$. And $(H,f)$  
is co-full at $v$ if and only if $\Phi_v^*(H,f)$ is 
full at $v$.   

In order to prove a duality theorem, we need to 
prepare a lemma.

\begin{lemma} Let $H$ and $K$ be Hilbert spaces and 
$T:H \rightarrow K$ be a bounded operator. Let 
$T= U|T|$ be its polar decomposition and 
$U$ a partial isometry 
with ${\it supp} \ U = \overline{\Im |T|}$ and 
$\Im U = \overline{\Im T}$. 
Suppose that 
$\Im T$ is closed. Then we have the following: 
\begin{itemize}
\item [(1)] 
$\Im |T| = \Im T^*$ is a closed subspace of $H$. 
\item [(2)]
Under the orthogonal decomposition 
\[
H =  \Ker |T|^{\perp} \oplus \Ker |T| 
  =  \Im |T|  \oplus \Ker |T|, 
\]
the restriction $|T||_{\Im |T|} : \Im |T| \rightarrow \Im |T|$ 
is a bounded invertible operator. 
\item[(3)] 
Let 
$S = (|T||_{\Im |T|})^{-1}$ be its inverse. 
Define a bounded operator $B: K \rightarrow \Im T^*$ by 
$Bx = SU^*x$ for $x \in K$. Let $Q :H \rightarrow \Im T^*$ 
be the canonical projection. 
Then $BT= Q$. Moreover $B|_{\Im T}:\Im T \rightarrow \Im T^*$ 
is a  bounded invertible operator.
\end{itemize}
\label{lemma:making-projection} 
\end{lemma}
\begin{proof}
(1)Since $\Im T$ is closed, $\Im T^*$ is also closed.  
Since $U(|T|x) = Tx$ by definition of $U$ and $\Im T$ is closed,  
 $\Im |T|$ is closed.  \\
(2)Since  $\Ker |T|^{\perp} = \Im |T|$, $|T||_{\Im |T|}$ is 
one to one. Since $|T|(H) = |T|(\Im |T|)$ is closed, 
$|T||_{\Im |T|}$ is onto. Hence $|T||_{\Im |T|}$ is 
bounded invertible. \\ 
(3)For any $x = x_1 + x_2 \in H$ 
with $x_1 \in \Im |T| = \Im T^* $ and $x_2 \in \Ker |T|$, 
\[
BTx = SU^*U|T|x = S|T|x = S|T|x_1 = x_1 = Qx . 
\]
It is clear that $B|_{\Im T}$ 
is a bounded invertible operator.
\end{proof}

\begin{thm}Let $\Gamma=(V,E,s,r)$ be a finite quiver and 
 $v \in V$ a sink. Assume that a Hilbert representation $(H,f)$ of 
$\Gamma$ is closed at $v$. 
Let $
h_v: \oplus_{\alpha \in E^v} H_{s(\alpha)} \rightarrow H_v 
$
be a bounded operator defined by 
$
h_v((x_{\alpha})_{\alpha \in E^v}) 
= \sum_{\alpha \in E^v} f_{\alpha}(x_{\alpha}).
$
Define  
a Hilbert representation $(\tilde{H},\tilde{f})$ of $\Gamma$
by 
$\tilde{H}_v = (\Im h_v)^{\perp} \subset H_v$,
$\tilde{H}_u = 0$ for $u \not=v$ and 
$\tilde{f} = 0$. 
Then we have 
\[
(H,f) \cong \Phi_v^{-}(\Phi_v^{+}(H,f)) \oplus (\tilde{H},\tilde{f}). 
\]
\label{thm:decomposition}
\end{thm} 
\begin{proof} 
Let $(H^+,f^+) = \Phi_v^{+}(H,f)$  and 
$(H^{+-},f^{+-}) = \Phi_v^{-}(\Phi_v^{+}(H,f))$. 
Then 
$H^+_v = \Ker h_v = \{ (x_{\alpha})_{\alpha \in E^v} \in 
\oplus_{\alpha \in E^v} H_{s(\alpha)}
 \ | \ 
      \sum_{\alpha \in E^v} f_{\alpha}(x_{\alpha}) = 0 \},$ 
and $H^+_u = H_u$ for $u \not=v$. 
We have  
$f^+_{\overline{\beta}}((x_{\alpha})_{\alpha \in E^v})
 = x_{\beta}$ for $\beta \in E^v$ 
, and $f^+_{\beta} = f_{\beta}$ for 
$\beta \not\in E^v$. 

Let $\hat{h}_v : H^+_v \rightarrow \oplus_{\alpha \in E^v} H_{s(\alpha)}
$ be a bounded operator given by 
\[
\hat{h}_v((x_{\alpha})_{\alpha \in E^v}) 
= (f^+_{\overline{\beta}}((x_{\alpha})_{\alpha \in E^v}))_{\beta \in E^v}
= (x_{\beta})_{\beta \in E^v}
= (x_{\alpha})_{\alpha \in E^v}. 
\]
Hence $\hat{h}_v$ is the canonical embedding. 
Since $(H,f)$ is closed at $v$, $\Im h_v$ and $\Im h_v^*$ 
are closed subspaces. 
Therefore 
\[
H^{+-}_v = (\Im \hat{h}_v)^{\perp} 
= (H^+_v)^{\perp}
= (\Ker h_v)^{\perp}
= \Im h_v^*. 
\]
For any other $u \in V$ with $u \not=v$, 
$H^{+-}_u = H_u$.  
Let $Q_v : \oplus_{\alpha \in E^v} H_{s(\alpha)}
 \rightarrow H^{+-}_v$ 
be the canonical projection.  
For $\beta \in E^v$, let 
\[
j_{\beta} : H_{s(\beta)} 
\rightarrow \oplus_{\alpha \in E^v} H_{s(\alpha)}
\]
be the canonical inclusion. 
Then 
$f^{+-}_{\beta} : H_{s(\beta)} \rightarrow H^{+-}_v$ 
is given by   $f^{+-}_{\beta} = Q_v \circ  j_{\beta} .  $ 
For other $\beta \not\in  E^v$, we have  
$f^{+-}_{\beta} = f_{\beta}$.
  
We shall define an isomorphism 
\[
\varphi : (H,f) \rightarrow 
\Phi_v^{-}(\Phi_v^{+}(H,f)) \oplus (\tilde{H},\tilde{f}).
\]
Apply Lemma \ref{lemma:making-projection}  by 
 putting  $T = h_v$, 
$H = \oplus_{\alpha \in E^v} H_{s(\alpha)}$ 
and $K = H_v$. Consider the polar decomposition 
$h_v = U|h_v|$. Put 
$S = (|h_v||_{\Im |h_v|})^{-1}$. 
Define a bounded operator $B: H_v\rightarrow \Im h_v^*$ by 
$B = SU^*$. Then $Bh_v$ 
is  the canonical projection $Q_v$ of  $H_v$ onto $\Im h_v^*$. 
We define 
\[
\varphi_v : H_v = \Im h_v \oplus (\Im h_v)^{\perp}
\rightarrow H^{+-}_v \oplus \tilde{H}_v 
=  \Im h_v^* \oplus (\Im h_v)^{\perp}
\]
by $\varphi_v(x,y) = (B|_{\Im h_v} x,y)$ for 
$x \in \Im h_v $ and 
$y \in (\Im h_v)^{\perp}$. 
By  Lemma \ref{lemma:making-projection} (2), 
$\varphi_v$ is a bounded invertible operator.  
For  $u \in V$ with $u \not=v$, put 
$\varphi _u : H_u \rightarrow H_u \oplus 0$ by 
$\varphi _u(x) = (x,0)$ for $x \in H_u$. 
For any $\beta\in E^v$ and $x \in H_{s(\beta)}$, 
\[
\varphi _v \circ f_{\beta}(x) 
=\varphi _v (h_v(j_{\beta}(x))) 
= (B(h_v(j_{\beta}(x))),0) 
= (Q_v(j_{\beta}(x)), 0). 
\]
On the other hand, 
\[
(f^{+-}_{\beta} \oplus 0) \circ \varphi _{s(\beta)} (x) 
= (f^{+-}_{\beta} \oplus 0)(x,0)
= (f^{+-}_{\beta}(x),0) 
= (Q_v \circ  j_{\beta}(x),0) .
\]
For other $\beta \not\in E^v$, we have 
\[
\varphi _{r(\beta)}\circ f^{+-}_{\beta} 
= \varphi _{r(\beta)}\circ f_{\beta} 
= f_{\beta} \circ \varphi _{s(\beta)} 
= f^{+-}_{\beta} \circ \varphi _{s(\beta)} .
\]
Hence 
$\varphi:(H,f) \rightarrow 
\Phi_v^{-}(\Phi_v^{+}(H,f)) \oplus (\tilde{H},\tilde{f})
$ is an isomorphism. 
\end{proof}

\smallskip
\noindent
{\bf Counter example. } If we do not assume that a Hilbert representation 
$(H,f)$ of $\Gamma$ is closed at $v$, then the above Theorem 
\ref{thm:decomposition} does not hold in general. In fact, 
consider the following quiver $\Gamma=(V,E,s,r)$:
\[
\circ_1 \overset{\alpha_1}\longrightarrow 
\circ_0 \overset{\alpha_2}\longleftarrow \circ_2 
\]
Let $K = \ell^2(\mathbb N)$ with the canonical basis 
$(e_n)_{n \in \mathbb N}$.  
Define a Hilbert representation $(H,f)$ of $\Gamma$ by 
$H_0= K \oplus K$, $H_1 = K \oplus 0$ and $H_2$ is the 
closed subspace  of $H_0$ spanned by 
$\{(\cos \frac{\pi}{n}e_n, \sin \frac{\pi}{n}e_n) \in K \oplus K 
\ | \ n \in \mathbb N\}$. Then $H_1 \cap H_2 = 0$ and 
$H_1 + H_2$ is a dense subspace of $H_0$ but not closed 
in $H_0$. Let 
$f_k = f_{\alpha_k}  : H_k \rightarrow  H_{0}$ be the 
inclusion map for $k = 1,2$. Then $(H,f)$ is not closed at 
a sink $0$. It is easy to see that $H^+_0 = \Ker h_0 = 0$, 
$f^+_1 = 0$ and $f^+_2 = 0$. Therefore 
$H^{+-}_0 = H_1 \oplus H_2$ and $H^{+-}_1 = H_1$, $H^{+-}_2 = H_2$.
We have $f^{+-}_k: H_k \rightarrow H_1 \oplus H_2$ is a canonical 
inclusion for $k = 1,2$. 
Since  $\tilde{H}_0 = (\Im h_v)^{\perp} = 0$, we have 
$(\tilde{H},\tilde{f}) = (0,0)$. Therefore 
\[ 
\Phi_v^{-}(\Phi_v^{+}(H,f)) \oplus (\tilde{H},\tilde{f})
= \Phi_v^{-}(\Phi_v^{+}(H,f)) 
= (H^{+-}, f^{+-})
\]
is closed at a sink $0$. But $(H,f)$ is not closed at a sink $0$. 
Therefore there exists no isomorphism between $(H,f)$ and 
$\Phi_v^{-}(\Phi_v^{+}(H,f)) \oplus (\tilde{H},\tilde{f})$. 

Note that $(H,f)$ is not full at a sink $0$ and 
$\Phi_v^{-}(\Phi_v^{+}(H,f))$ is full at a sink $0$. Therefore 
this example also shows that,  
if we do not assume that a Hilbert representation 
$(H,f)$ of $\Gamma$ is full at $v$, then the following 
Duality Theorem (Corollary 
\ref{cor:duality}) does not hold in general.

\begin{cor}{\rm (Duality theorem.)}
Let $\Gamma=(V,E,s,r)$ be a finite quiver and 
 $v \in V$ a sink. If  a Hilbert representation $(H,f)$ of 
$\Gamma$ is full at $v$, then 
\[
(H,f) \cong \Phi_v^{-}(\Phi_v^{+}(H,f)). 
\]
\label{cor:duality}
\end{cor}
\begin{proof}
Since $(H,f)$ is full at $v$, 
$\tilde{H}_v = (\Im h_v)^{\perp}= H_v^{\perp} = 0$. Hence 
$(\tilde{H},\tilde{f}) = (0,0)$ in 
Theorem \ref{thm:decomposition}. 
\end{proof}

\noindent
{\bf Remark.} It is also necessary that 
$(H,f)$ is full at the sink $v$ in order that the 
above Duality Theorem 
holds. It follows from Lemma \ref{lemma:full} below. 

We have a dual version. 

\begin{thm}Let $\Gamma=(V,E,s,r)$ be a finite quiver and 
 $v \in V$ a source. Assume that a Hilbert representation $(H,f)$ of 
$\Gamma$ is co-closed at $v$. 
Let 
$
\hat{h}_v: H_v  \rightarrow \oplus_{\alpha \in E_v} H_{r(\alpha)}   
$
is a bounded operator defined by 
$
\hat{h}_v(x) 
= (f_{\alpha}(x))_{\alpha \in E_v} \ \text{ for } x \in H_v.
$
Define 
a Hilbert representation $(\check{H},\check{f})$ of $\Gamma$
by 

\[
\check{H}_v = (\Im \hat{h}_v^*)^{\perp} 
(= \Ker \hat{h}_v = 
\cap_{\alpha \in E_v} \Ker f_{\alpha} ) \subset H_v,
\] 
$\tilde{H}_u = 0$ for $u \not=v$ and 
$\tilde{f} = 0$. 
Then 
\[
(H,f) \cong \Phi_v^{+}(\Phi_v^{-}(H,f)) \oplus (\check{H},\check{f}). 
\]
\label{thm:decomposition2}
\end{thm} 
\begin{proof} We see that $v$ is a sink in  $\overline{\Gamma}$,
because $v$ is a source in $\Gamma$. 
Since a Hilbert representation $(H,f)$ of 
$\Gamma$ is co-closed at $v$, a Hilbert representation
$\Phi_v^*(H,f)$ is 
closed at $v$. By Theorem \ref{thm:decomposition}, 
there exists  
a Hilbert representation $(\tilde{H},\tilde{f})$ 
of $\overline{\Gamma}$
such that 
\[
\Phi_v^*(H,f) \cong \Phi_v^{-}(\Phi_v^{+}(\Phi_v^*(H,f))) 
\oplus (\tilde{H},\tilde{f}). 
\]
Put $(\check{H},\check{f})= \Phi_v^*(\tilde{H},\tilde{f})$. Then
\begin{align*}
& (H,f) \cong 
\Phi_v^*(\Phi_v^*(H,f))
\cong \Phi_v^*\Phi_v^{-}\Phi_v^{+}\Phi_v^*(H,f)
\oplus \Phi_v^*(\tilde{H},\tilde{f}) \\
& \cong (\Phi_v^*\Phi_v^{-}\Phi_v^*)(\Phi_v^*\Phi_v^{+}\Phi_v^*)(H,f) 
\oplus \Phi_v^*(\tilde{H},\tilde{f})   \\
& \cong 
\Phi_v^{+}(\Phi_v^{-}(H,f)) \oplus (\check{H},\check{f}).
\end{align*}
Moreover it is easy to see that 
\[
\check{H}_v = ( \sum_{\alpha \in E_v} \Im f_{\alpha}^*   )^{\perp}
= \cap_{\alpha \in E_v} \Ker f_{\alpha}. 
\]
\end{proof}

\smallskip
\noindent
{\bf Counter example. } If we do not assume that a Hilbert representation 
$(H,f)$ of $\Gamma$ is co-closed at the source 
$v$, then the above Theorem 
\ref{thm:decomposition2} does not hold in general. In fact, 
consider the following quiver $\Gamma=(V,E,s,r)$:
\[
\circ_1 \overset{\alpha_1}\longleftarrow 
\circ_0 \overset{\alpha_2}\longrightarrow \circ_2 
\]
Let $K = \ell^2(\mathbb N)$ with the canonical basis 
$(e_n)_{n \in \mathbb N}$.  
Define a Hilbert representation $(H,f)$ of $\Gamma$ by 
$H_0= K \oplus K$, $H_1 = K \oplus 0$ and $H_2$ is the 
closed subspace  $H_0$ spanned by 
$\{(\cos \frac{\pi}{n}e_n, \sin \frac{\pi}{n}e_n) \in K \oplus K 
\ | \ n \in \mathbb N\}$.  Let 
$f_k = f_{\alpha_k}  : H_0 \rightarrow  H_k$ be the canonical 
projection for $k = 1,2$. Then $(H,f)$ is not co-closed at 
a source $0$. It is easy to see that 
$H^-_0 = (\Im {\hat h}_0)^{\perp} = 0$, 
$f^-_1 = 0$ and $f^-_2 = 0$. Therefore 
$H^{-+} _0 = H_1 \oplus H_2$ and $H_1^{-+}  = H_1$, $H_2^{-+} = H_2$. 
We have that $f_k^{-+} : H_1 \oplus H_2 \rightarrow H_k$ 
is the canonical 
projection for $k = 1,2$. 
Since  $\check{H_0} = \Ker {\hat h}_0 = 0$, we have 
$(\check{H},\check{f}) = (0,0)$. Therefore 
\[ 
\Phi_v^{+}(\Phi_v^{-}(H,f)) \oplus (\check{H},\check{f})
= \Phi_v^{+}(\Phi_v^{-}(H,f)) 
=(H^{-+}, f^{-+})
\]
is co-closed at a source $0$. But $(H,f)$ is not co-closed 
at a source $0$. 
Therefore there exists no isomorphism between $(H,f)$ and 
$\Phi_v^{+}(\Phi_v^{-}(H,f)) \oplus (\check{H},\check{f})$. 

Note that $(H,f)$ is not co-full at a source $0$ and 
$\Phi_v^{+}(\Phi_v^{-}(H,f))$is co-full at a source $0$. Therefore 
this example also shows that,  
if we do not assume that a Hilbert representation 
$(H,f)$ of $\Gamma$ is co-full at $v$, then the following 
Duality Theorem (Corollary 
\ref{cor:duality2}) does not hold in general.

\begin{cor}{\rm (Duality theorem.)}
Let $\Gamma=(V,E,s,r)$ be a finite quiver and 
 $v \in V$ a source. If  a Hilbert representation $(H,f)$ of 
$\Gamma$ is co-full at $v$, then 
\[
(H,f) \cong \Phi_v^{+}(\Phi_v^{-}(H,f)). 
\]
\label{cor:duality2}
\end{cor}
\begin{proof}
Since $(H,f)$ is co-full at $v$, 
$\check{H}_v = \cap_{\alpha \in E_v} \Ker f_{\alpha} = 0$. Hence 
$(\check{H},\check{f}) = (0,0)$ in 
Theorem \ref{thm:decomposition2}. 
\end{proof}

\noindent
{\bf Remark.} It is also necessary that 
$(H,f)$ is co-full at the source $v$ in order that the 
above Duality Theorem 
holds. It follows from Lemma \ref{lemma:co-full} below.

\begin{lemma}
Let $\Gamma=(V,E,s,r)$ be a finite quiver and 
 $v \in V$ a sink. Then for any  Hilbert representation $(H,f)$ of 
$\Gamma$,  $\Phi_v^{+}(H,f)$ is co-full at $v$. 
\label{lemma:co-full}
\end{lemma}
\begin{proof} Put $(H^+,f^+) = \Phi_v^{+}(H,f)$. 
Recall that 
$h_v: \oplus_{\alpha \in E^v} H_{s(\alpha)} \rightarrow H_v $
is given by 
$h_v((x_{\alpha})_{\alpha \in E^v}) 
= \sum_{\alpha \in E^v} f_{\alpha}(x_{\alpha})$, 
and $H^+_v = \Ker h_v$. And 
For $\beta \in E^v$, 
let  $i_v : H^+_v \rightarrow \oplus_{\alpha \in E^v} H_{s(\alpha)}$
be the canonical inclusion
and $P_{\beta} : \oplus_{\alpha \in E^v} H_{s(\alpha)} 
\rightarrow H_{s(\beta)}$ 
the canonical projection. We define 
\[
f^+_{\overline{\beta}} :  H^+_{s(\overline{\beta})} = H^+_v \rightarrow 
H^+_{r(\overline{\beta})} = H_{s(\beta)} \ \ 
\text{ by } \  g_{\overline{\beta}} = P_{\beta} \circ i_v. 
\]
Therefore 
${f^+_{\overline{\beta}}}^* :  H_{s(\beta)}  \rightarrow H^+_v$ 
is given by ${f^+_{\overline{\beta}}}^* = i_v^* \circ {P_{\beta}}^*$. 
Since  
$P_{\beta}^* : H_{s(\beta)} \rightarrow 
\oplus_{\alpha \in E^v} H_{s(\alpha)}$ 
is the canonical inclusion and 
$i_v^*: \oplus_{\alpha \in E^v} H_{s(\alpha)} 
\rightarrow H^+_v$ is the canonical projection, 
we have 
\[
\sum_{\overline{\beta} \in E_v}   \Im {f^+_{\overline{\beta}}}^* 
 = \sum_{\beta \in E^v} \Im (i_v^* \circ {P_{\beta}}^*) 
 = H^+_v . 
\]
Therefore $(H^+,f^+)$ is co-full at $v$. 
\end{proof}

\begin{prop}
Let $\Gamma=(V,E,s,r)$ be a finite quiver and 
 $v \in V$ a sink. If $(H,f)$ is a Hilbert representation  of 
$\Gamma$, then 
\[
 \Phi_v^{+} \Phi_v^{-}\Phi_v^{+}(H,f) \cong \Phi_v^{+}(H,f). 
\]
\end{prop}
\begin{proof}
Since $\Phi_v^{+}(H,f)$ is co-full at the source $v$ 
in $\sigma_v^{+}(\Gamma)$ by the above lemma \ref{lemma:co-full}, 
duality theorem (Corollary 
\ref{cor:duality2} ) yields the conclusion.  
\end{proof}

\begin{lemma}
Let $\Gamma=(V,E,s,r)$ be a finite quiver and 
 $v \in V$ a source. Then for any  Hilbert representation $(H,f)$ of 
$\Gamma$,  $\Phi_v^{-}(H,f)$ is full at $v$. 
\label{lemma:full}
\end{lemma}
\begin{proof} Put $(H^-,f^-) = \Phi_v^{-}(H,f)$. 
Recall that 
$
\hat{h}_v: H_v  \rightarrow \oplus_{\alpha \in E_v} H_{r(\alpha)}   
$
is given  by 
$
\hat{h}_v(x) 
= (f_{\alpha}(x))_{\alpha \in E_v} \ \text{ for } x \in H_v.
$
and 
$H^-_v =  (\Im \hat{h}_v)^{\perp}  \subset 
\oplus_{\alpha \in E_v} H_{r(\alpha)},  
$
Let $Q_v : \oplus_{\alpha \in E_v} H_{r(\alpha)} \rightarrow H^-_v$ 
be the canonical projection.  
For $\beta \in E_v$, let 
$
j_{\beta} : H_{r(\beta)} 
\rightarrow \oplus_{\alpha \in E_v} H_{r(\alpha)}
$
be the canonical inclusion. 
Then 
\[
f^-_{\overline{\beta}} :  H^-_{s(\overline{\beta})} = 
H_{r(\beta)} \rightarrow H^-_{r(\overline{\beta})} = H_v ^- \ \ 
\text{ by }  f^-_{\overline{\beta}} = Q_v \circ j_{\beta} .
\]
Therefore 
\[
\sum_{\overline{\beta} \in E^v } \Im f^-_{\overline{\beta}} 
= Q_v(\oplus_{\alpha \in E_v} H_{r(\alpha)}) 
    = H^-_v . 
\]

Thus  $(H^-,f^-)$ is full at $v$. 
\end{proof}

\begin{prop}
Let $\Gamma=(V,E,s,r)$ be a finite quiver and 
 $v \in V$ a source. If $(H,f)$ is a Hilbert representation  of 
$\Gamma$, then 
\[
 \Phi_v^{-} \Phi_v^{+}\Phi_v^{-}(H,f) \cong \Phi_v^{-}(H,f). 
\]
\end{prop}
\begin{proof}
Since $\Phi_v^{-}(H,f)$ is full at the source 
in $\sigma_v^{-} (\Gamma)$ by the above lemma \ref{lemma:full}, 
duality theorem (Corollary 
\ref{cor:duality} )
 yields the conclusion.  
\end{proof}

We examine on which representation a reflection functor vanishes. 

\begin{lemma}
Let $\Gamma=(V,E,s,r)$ be a finite quiver and 
 $v \in V$ a sink. Then,  for any  Hilbert representation $(H,f)$ of 
$\Gamma$, the following are equivalent: 
\begin{itemize}
\item [(1)] 
$\Phi_v^{+}(H,f) =\cong (0,0)$ 
\item[(2)]
$H_u = 0$ for any $u \in V$ with $u \not= v$.  
\end{itemize}
Furthermore if the above conditions are satisfied and 
$(H,f)$ is indecomposable, then $H_v \cong \mathbb C$.  
If the above conditions are satisfied and 
$(H,f)$ is full at the sink $v$, then $(H,f) \cong (0,0)$.  

\label{lemma:vanish}
\end{lemma}
\begin{proof}
Put $(H^+,f^+) = \Phi_v^{+}(H,f)$. 
Recall that 
$h_v: \oplus_{\alpha \in E^v} H_{s(\alpha)} \rightarrow H_v $ 
is given by 
$h_v((x_{\alpha})_{\alpha \in E^v}) 
= \sum_{\alpha \in E^v} f_{\alpha}(x_{\alpha})$, 
and $H^+_v = \Ker h_v$. For other $u \in V$ with $u \not= v$, 
$H^+_u = H_u$. 

\noindent
(1)$\Rightarrow$(2):Assume that 
$\Phi_v^{+}(H,f) = 0$.  Then,  for any  $u \in V$ with $u \not= v$ 
we have $H_u = H^+_u = 0$. \\
\noindent
(2)$\Rightarrow$(1): Assume that  
$H_u = 0$ for any $u \in V$ with $u \not= v$. 
Then $H^+_v = 0$, because $H^+_v 
= \Ker h_v \subset \oplus_{\alpha \in E^v} H_{s(\alpha)} = 0$. 
 For other $u \in V$ with $u \not= v$, 
$H^+_u = H_u = 0$. 

Furthermore assume that the above conditions are satisfied 
and $(H,f)$ is indecomposable.  Then $f = 0$. 
Suppose that  $\dim H_v \geq 2$. Then a non-trivial decomposition 
$H_v = K \oplus L$ gives a non-trivial decomposition of $(H,f)$. 
This contradicts that $(H,f)$ is indecomposable. Hence 
$H_v \cong \mathbb C$.
Assume that  the above conditions are satisfied 
and $(H,f)$ is full at $v$.  Then $f = 0$, so that 
$H_v = \sum_{\alpha \in E^v} \Im f_{\alpha} = 0$. 
Hence $(H,f) \cong (0,0)$.

\end{proof}

\begin{lemma}
Let $\Gamma=(V,E,s,r)$ be a finite quiver and 
 $v \in V$ a source. Then,  for any  Hilbert representation $(H,f)$ of 
$\Gamma$, the following condition  are equivalent: 
\begin{itemize}
\item [(1)] 
$\Phi_v^{-}(H,f) \cong (0,0)$ 
\item[(2)]
$H_u = 0$ for any $u \in V$ with $u \not= v$.  
\end{itemize}
Furthermore if the above conditions are satisfied and 
$(H,f)$ is indecomposable, then $H_v \cong \mathbb C$. 
If the above conditions are satisfied and 
$(H,f)$ is co-full at the source $v$ , then $(H,f) \cong (0,0)$.  
 
\label{lemma:vanish2}
\end{lemma}
\begin{proof}
Put $(H^-,f^-) = \Phi_v^{-}(H,f)$. 
Recall that 
$
\hat{h}_v: H_v  \rightarrow \oplus_{\alpha \in E_v} H_{r(\alpha)}   
$ is given  by 
$
\hat{h}_v(x) 
= (f_{\alpha}(x))_{\alpha \in E_v} \ \text{ for } x \in H_v, 
$
and 
$H^-_v =  (\Im \hat{h}_v)^{\perp}  \subset 
\oplus_{\alpha \in E_v} H_{r(\alpha)}$.   
For other $u \in V$ with $u \not= v$, 
$H^-_u = H_u$. 

\noindent
(1)$\Rightarrow$(2):Assume that 
$\Phi_v^{-}(H,f) = 0$.  Then,  for any  $u \in V$ with $u \not= v$ 
we have $H_u = H^-_u = 0$. \\
\noindent
(2)$\Rightarrow$(1): Assume that  
$H_u = 0$ for any $u \in V$ with $u \not= v$. 
Then $H^-_v = 0$, because $H^-_v 
=  (\Im \hat{h}_v)^{\perp} \subset 
\oplus_{\alpha \in E_v} H_{r(\alpha)} = 0$.
For other $u \in V$ with $u \not= v$, 
$H^-_u = H_u = 0$. 

Assume that the above conditions are satisfied 
and $(H,f)$ is co-full at $v$.  Since $f_{\alpha}^* = 0$ 
for any $\alpha \in E$,   
$H_v = \sum_{\alpha \in E_v} \Im f_{\alpha}^* = 0$. 
Hence $(H,f) \cong (0,0)$.
The rest is clear. 
\end{proof}

We shall show that a reflection  functor preserves indecomposability 
of a Hilbert representation unless vanishing  on it, under the 
assumption that  the 
Hilbert representation is closed (resp. co-closed) 
at a sink (resp. source).   

\begin{thm}
Let $\Gamma=(V,E,s,r)$ be a finite quiver and 
 $v \in V$ a sink. Suppose that 
 a Hilbert representation $(H,f)$ of $\Gamma$ is indecomposable 
and  closed at $v$. Then we have the following: 
\begin{itemize}
\item [(1)] 
If $\Phi_v^{+}(H,f) = 0$, then $H_v = \mathbb C$, 
$H_u = 0$ for any $u \in V$ with $u \not= v$ and 
$f_{\alpha} = 0$ for any $\alpha \in E$. 
\item[(2)]
If $\Phi_v^{+}(H,f) \not= 0$, then $\Phi_v^{+}(H,f)$ 
is also indecomposable  and 
$(H,f) \cong \Phi_v^{-}(\Phi_v^{+}(H,f)). $
\label{thm:preserve}
\end{itemize}
\end{thm}

\begin{proof}
Recall an operator $
h_v: \oplus_{\alpha \in E^v} H_{s(\alpha)} \rightarrow H_v 
$
defined by 
$
h_v((x_{\alpha})_{\alpha \in E^v}) 
= \sum_{\alpha \in E^v} f_{\alpha}(x_{\alpha}).
$
Since $(H,f)$ is closed at a sink $v$, we have a 
decomposition 
such that 
\[
(H,f) \cong \Phi_v^{-}(\Phi_v^{+}(H,f)) \oplus (\tilde{H},\tilde{f})
\]
by Theorem \ref{thm:decomposition}, 
where 
$\tilde{H}_v = (\Im h_v)^{\perp} \subset H_v$,
$\tilde{H}_u = 0$ for $u \not=v$ and 
$\tilde{f} = 0$. 

Since $(H,f)$ is indecomposable,   
$\Phi_v^{-}(\Phi_v^{+}(H,f)) \cong (0,0)$ or 
$(\tilde{H},\tilde{f}) \cong (0,0)$. 

\noindent
(Case 1): Suppose that  $\Phi_v^{-}(\Phi_v^{+}(H,f)) \cong (0,0)$. 
Then $(H,f) \cong (\tilde{H},\tilde{f})$. Hence 
$H_u \cong  \tilde{H}_u = 0$ for $u \not=v$. This implies that 
$\Phi_v^{+}(H,f) \cong (0,0)$ by Lemma \ref{lemma:vanish}. 
Since $(H,f)$ is indecomposable, $H_v \cong \mathbb C$. 

\noindent
(Case 2):Suppose that $(\tilde{H},\tilde{f}) \cong (0,0)$. 
Then $(H,f) \cong \Phi_v^{-}(\Phi_v^{+}(H,f))$. Since 
$(H,f)$ is non-zero, $\Phi_v^{+}(H,f)$ is non-zero. We 
shall show that $\Phi_v^{+}(H,f)$ is indecomposable. 
Assume that $\Phi_v^{+}(H,f)\cong (K,g) \oplus (K',g')$. 
Then 
\[
(H,f) \cong \Phi_v^{-}(\Phi_v^{+}(H,f)) \cong 
 \Phi_v^{-}(K,g) \oplus \Phi_v^{-}(K',g') .
\]
Since $(H,f)$ is indecomposable, 
$\Phi_v^{-}(K,g) \cong (0,0)$ or $\Phi_v^{-}(K',g') \cong (0,0)$.
By Lemma \ref{lemma:co-full}
, $\Phi_v^{+}(H,f)$ is co-full at $v$, so are its direct 
summands $(K,g)$ and $(K',g')$. Then $(K,g) \cong (0,0)$ or 
$(K',g') \cong (0,0)$ by Lemma \ref{lemma:vanish2}. Thus 
$\Phi_v^{+}(H,f)$ is indecomposable. 

Since (Case 1) and (Case 2) are mutually exclusive and 
either of them occurs, we get the conclusion. 

\end{proof}

We have a dual version.

\begin{thm}
Let $\Gamma=(V,E,s,r)$ be a finite quiver and 
 $v \in V$ a source. Suppose that 
 a Hilbert representation $(H,f)$ of $\Gamma$ is indecomposable 
and  co-closed at $v$. Then we have the following: 
\begin{itemize}
\item [(1)] 
If $\Phi_v^{-}(H,f) = 0$, then $H_v = \mathbb C$, 
$H_u = 0$ for any $u \in V$ with $u \not= v$ and 
$f_{\alpha} = 0$ for any $\alpha \in E$. 
\item[(2)]
If $\Phi_v^{-}(H,f) \not= 0$, then $\Phi_v^{-}(H,f)$ 
is also indecomposable  and 
$(H,f) \cong \Phi_v^{+}\Phi_v^{-}(H,f)). $
\end{itemize}
\label{thm:preserving-indecomposability} 

\end{thm}
\begin{proof}
A dual argument of the proof in Theorem \ref{thm:preserve} works.  

\end{proof} 

\section{Extended Dynkin diagrams}

Gabriel's theorem says that a connected finite quiver has 
only finitely many indecomposable representations if and only if 
the underlying undirected graph is one of Dynkin diagrams 
$A_n, D_n, E_6, E_7,E_8$. Representations 
of quivers on finite-dimensional vector spaces has 
been developed by  Bernstein-Gelfand-Ponomarev \cite{BGP}, 
Donovan-Freislish \cite{DF}, V. Dlab-Ringel \cite{DR}, 
Gabriel-Roiter \cite{GR}, 
Kac \cite{Ka}, Nazarova \cite{Na} ... . 

Furthermore locally  scalar representations of quivers 
in the category of Hilbert spaces up to the unitary 
equivalence were introduced by 
Kruglyak and Roiter \cite{KRo}. They prove an 
analog of Gabriel's theorem.

We consider 
a complement of  Gabriel's theorem for Hilbert 
representations. 
We need to construct some examples of indecomposable, 
infinite-dimensional
representations of quivers  with the 
underlying  undirected graphs extended Dynkin diagrams 
$\tilde{D_n} \  (n \geq 4), ,\tilde{E_7}$ and $\tilde{E_8}$. 
We consider the relative position of several subspaces along 
the quivers, where 
vertices are represented by a family of subspaces and 
 arrows are  represented by natural inclusion maps.

\begin{lemma}
Let $\Gamma=(V,E,s,r)$ be the following quiver with the 
underlying  undirected graph an extended Dynkin diagram 
$\tilde{D_n}$ for $n \geq 4$:

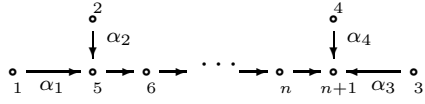
\begin{picture}(150,60)(-75,5)
\put(0,25){\thicklines\circle{2}}
\put(0,20){${}_{{}_{1}}$}

\put(5,25){\vector(1,0){20}}

\put(10,20){${}_{\alpha_{1}}$}

\put(30,25){\thicklines\circle{2}}

\put(30,20){${}_{{}_{5}}$}

\put(30,40){\vector(0,-1){10}}

\put(30,45){\thicklines\circle{2}}

\put(30,50){${}_{{}_{2}}$}

\put(35,38){${}_{\alpha_{2}}$}

\put(35,25){\vector(1,0){10}}

\put(50,25){\thicklines\circle{2}}

\put(50,20){${}_{{}_{6}}$}

\put(55,25){\vector(1,0){10}}

\put(70,25){$\cdots$}

\put(85,25){\vector(1,0){10}}

\put(100,25){\thicklines\circle{2}}

\put(100,20){${}_{{}_{n}}$}

\put(105,25){\vector(1,0){10}}

\put(120,25){\thicklines\circle{2}}

\put(115,20){${}_{{}_{n+1}}$}

\put(145,25){\vector(-1,0){20}}

\put(134,20){${}_{\alpha_{3}}$}

\put(150,25){\thicklines\circle{2}}

\put(150,20){${}_{{}_{3}}$}

\put(120,40){\vector(0,-1){10}}

\put(120,45){\thicklines\circle{2}}

\put(120,50){${}_{{}_{4}}$}

\put(125,38){${}_{\alpha_{4}}$}

\end{picture}

\noindent
Then there exists an infinite-dimensional, 
indecomposable Hilbert representation $(H,f)$ of $\Gamma$.  
\label{lemma:Dn}
\end{lemma}
\begin{proof}
Let $K = \ell^2(\mathbb N)$ and $S$ a unilateral shift on $K$. 
We define a Hilbert representation
$(H,f) := ((H_v)_{v\in V},(f_{\alpha})_{\alpha \in E})$ 
of $\Gamma$ as follows: \\
Define 
 $H_1 = K \oplus 0$, 
$H_2 = 0 \oplus K$, 
$H_3 = \{(x,Sx)\in K \oplus K | x \in K\}$, \\
$H_4 = \{(x,x)\in K \oplus K | x \in K\}$. 
$H_5 = H_6 = \dots =H_{n+1} = K \oplus K$, \\
Let $f_{\alpha_k}  
: H_{s(\alpha_k)}  \rightarrow  H_{r(\alpha_k)} $ be 
the inclusion map for any  $\alpha_k \in E$ for 
$k = 1,2,3,4$,  and $f_{\beta} = id$ for other arrows 
$\beta \in E$. Then 
we can show that $(H,f)$ is indecomposable
as in Example 3 in section 3. 

\end{proof}

Let $\Gamma=(V,E,s,r)$ be the quiver of  Example 4 in section 3.
with  the 
underlying  undirected graph a extended Dynkin diagram 
$\tilde{E_6}$. 
We have already shown that 
there exists an infinite-dimensional, 
indecomposable Hilbert representation $(H,f)$ of $\Gamma$. 

\begin{lemma}
Let $\Gamma=(V,E,s,r)$ be the following quiver with the 
underlying  undirected graph an extended Dynkin diagram 
$\tilde{E_7}$:

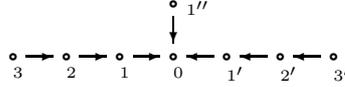
\begin{picture}(140,45)(-80,0)

\put(10,15){\thicklines\circle{2}}

\put(10,10){${}_{{}_{3}}$}

\put(15,15){\vector(1,0){10}}

\put(30,15){\thicklines\circle{2}}

\put(30,10){${}_{{}_{2}}$}

\put(35,15){\vector(1,0){10}}

\put(50,15){\thicklines\circle{2}}

\put(50,10){${}_{{}_{1}}$}

\put(55,15){\vector(1,0){10}}

\put(70,15){\thicklines\circle{2}}

\put(70,10){${}_{{}_{0}}$}

\put(85,15){\vector(-1,0){10}}

\put(90,15){\thicklines\circle{2}}

\put(90,10){${}_{{}_{1^{\prime}}   }$}

\put(105,15){\vector(-1,0){10}}

\put(110,15){\thicklines\circle{2}}

\put(110,10){${}_{{}_{2^{\prime} }  }$}

\put(125,15){\vector(-1,0){10}}

\put(130,15){\thicklines\circle{2}}

\put(130,10){${}_{{}_{3^{\prime} }  }$}

\put(70,30){\vector(0,-1){10}}

\put(70,35){\thicklines\circle{2}}

\put(75,35){${}_{{}_{ 1^{{\prime}{\prime}} }  }$}

\end{picture}

\noindent
Then there exists an infinite-dimensional, 
indecomposable Hilbert representation $(H,f)$ of $\Gamma$.  
\label{lemma:E7}
\end{lemma}
\begin{proof}
Let $K = \ell^2(\mathbb N)$ and $S$ a unilateral shift on $K$. 
We define a Hilbert representation
$(H,f) := ((H_v)_{v\in V},(f_{\alpha})_{\alpha \in E})$ 
of $\Gamma$ as follows: \\
Let 
$H_{0}=K\oplus K\oplus K\oplus K,$
$H_{1}=K\oplus0\oplus K\oplus K,$ \\
$H_{2}=K\oplus0\oplus\{(x,x);x\in K\},$
$H_{3}=K\oplus0\oplus0\oplus0,$ \\
$H_{1^{^{\prime}}}=0\oplus K\oplus K\oplus K,$
$H_{2^{^{\prime}}}=0\oplus K\oplus\{(y,Sy)\in K^2 \ | y\in K\},$ \\
$H_{3^{^{\prime}}}=0\oplus K\oplus0\oplus0$ and 
$H_{1^{^{\prime\prime}}}=\{(x,y,x,y)\in K^4 \ | \ x,y \in K\}.$
For any arrow $\alpha \in E$, let  
$f_{\alpha} : H_{s(\alpha)} \rightarrow H_{r(\alpha)}$ be 
the canonical inclusion map. We shall show that 
the Hilbert representation $(H,f)$ is indecomposable. 
Take $T = (T_v)_{v \in V} \in Idem (H,f)$. 
Since $T \in End (H,f)$ and any arrow is represented 
by the inclusion map, we have $T_0 x = T_vx$ for any 
$v \in \{1,2,3,1',2',3',1'' \}$ and 
any  $x \in H_v$. In particular, 
$T_0 H_v \subset H_v$.  
Since $T_{0}$ preserves $H_{3}=K\oplus0\oplus0\oplus0$, 
$H_{3^{^{\prime}}}=0\oplus K\oplus0\oplus0,$ 
and $H_{1^{^{\prime}}}\cap H_{1}=0\oplus0\oplus K\oplus K$, 
 $T_{0}$ is written 
\[
T_{0}=\left(
\begin{array}
[c]{cccc}%
A & 0 & 0 & 0\\
0 & B & 0 & 0\\
0 & 0 & X & Y\\
0 & 0 & Z & W
\end{array}
\right)  , 
\]
for some $A,B,X,T,Z,W \in B(K)$. 

Because  $H_{1^{^{\prime\prime}}}
=\{(x,y,x,y)\in K^4 \ | \ x,y\in K\}$ 
is also invariant under $T_0$, for any $x,y \in K$, 
there exist $x', y' \in K$ such that 
\[
\left(
\begin{array}
[c]{cccc}%
A & 0 & 0 & 0\\
0 & B & 0 & 0\\
0 & 0 & X & Y\\
0 & 0 & Z & W
\end{array}
\right)  \left(
\begin{array}
[c]{c}%
x\\
y\\
x\\
y
\end{array}
\right)  =\left(
\begin{array}
[c]{c}%
Ax\\
By\\
Xx+Yy\\
Zx+Wy
\end{array}
\right)  =\left(
\begin{array}
[c]{c}%
x^{^{\prime}}\\
y^{^{\prime}}\\
x^{^{\prime}}\\
y^{^{\prime}}%
\end{array}
\right)  .
\]
Putting $y = 0$, 
we have $Ax = Xx$ and $0 = Zx$ for any $x \in K$. 
Hence $A = X$ and $Z =0$. Similarly,  letting $x = 0$, 
we have $Y = 0$ and $W = B$.  
Therefore $T_0$ has a block diagonal form such that 
\[
T_{0}=\left(
\begin{array}
[c]{cccc}%
A & 0 & 0 & 0\\
0 & B & 0 & 0\\
0 & 0 & A & 0\\
0 & 0 & 0 & B
\end{array}
\right)  =A\oplus B\oplus A\oplus B.
\]
Furthermore, as  $T_0$ preserves 
$H_{1'} \cap H_{2}= \{(0,0,x,x)\in K^4 \ | \ x\in K\}, $
for any $x \in K$ there exists $y \in K$ such that 
$(0,0,Ax,Bx) = (0,0,y,y)$. 
Hence $A = B$. Therefore 
$T_{0}=A\oplus A\oplus A\oplus A.$
Moreover 
$H_1 \cap H_{2^{^{\prime}}}
= \{(0,0,x,Sx) \in K^4 \ | \ x \in K \}$
is also invariant under $T_0$. 
Hence for any $x \in K$, there exists $y \in K$ such that 
$(0,0,Ax,ASx) = (0,0,y,Sy)$. Thus $AS=SA$. 
Since $T\in Idem(H,f)$, 
$T_{0}$ is idempotent, so that $A$ is also idempotent.
Because $AS=SA$ and $A^2 = A$, we have 
$A = 0$ or $A = I$. Thus $T_0 = 0$ or 
$T_0 = I$. 
Since for any $v \in V$ and any  $x \in H_v$ 
$T_0 x = T_vx$, 
we have $T_v=0$ or $T_v=I$ simultaneously. 
Thus $T = (T_v)_{v \in V} = 0$ or $T = I$, that is, 
$Idem(H,f)=\{0,I\}.$
Therefore $(H,f)$ is indecomposable.
\end{proof}

\noindent
{\bf Remark.} Replacing $S$ by $S + \lambda I$ for $\lambda 
\in {\mathbb C}$, 
 we have uncountably  many 
 infinite-dimensional, indecomposable  
 Hilbert representations of  $\Gamma$.

\begin{lemma}
Let $\Gamma=(V,E,s,r)$ be the following quiver with the 
underlying  undirected graph an extended Dynkin diagram 
$\tilde{E_8}$:

\begin{picture}(150,45)(-70,0)

\put(10,15){\thicklines\circle{2}}

\put(10,10){${}_{{}_{5}}$}

\put(15,15){\vector(1,0){10}}

\put(30,15){\thicklines\circle{2}}

\put(30,10){${}_{{}_{4}}$}

\put(35,15){\vector(1,0){10}}

\put(50,15){\thicklines\circle{2}}

\put(50,10){${}_{{}_{3}}$}

\put(55,15){\vector(1,0){10}}

\put(70,15){\thicklines\circle{2}}

\put(70,10){${}_{{}_{2}}$}

\put(75,15){\vector(1,0){10}}

\put(90,15){\thicklines\circle{2}}

\put(90,10){${}_{{}_{1} }$}

\put(95,15){\vector(1,0){10}}

\put(110,15){\thicklines\circle{2}}

\put(110,10){${}_{{}_{0} }$}

\put(125,15){\vector(-1,0){10}}

\put(130,15){\thicklines\circle{2}}

\put(130,10){${}_{{}_{1^{\prime} }  }$}

\put(145,15){\vector(-1,0){10}}

\put(150,15){\thicklines\circle{2}}

\put(150,10){${}_{{}_{2^{\prime} }  }$}

\put(110,35){\thicklines\circle{2}}
\put(115,35){${}_{{}_{ 1^{{\prime}{\prime}} }  }$}
\put(110,30){\vector(0,-1){10}}

\end{picture}

\noindent
Then there exists an infinite-dimensional, 
indecomposable Hilbert representation $(H,f)$ of $\Gamma$.  
\label{lemma:E8}
\end{lemma}

\begin{proof}
Let $K = \ell^2(\mathbb N)$ and $S$ a unilateral shift on $K$. 
We define a Hilbert representation
$(H,f) := ((H_v)_{v\in V},(f_{\alpha})_{\alpha \in E})$ 
of $\Gamma$ as follows: \\
Let 
$H_{0}=K\oplus K\oplus K\oplus K\oplus K\oplus K,$ \\
$H_{1}=\{(x,x)\in K^2 \ | \ x\in K\} 
\oplus K\oplus K\oplus K\oplus K,$ \\
$H_{2}=0\oplus0\oplus K\oplus K\oplus K\oplus K,$
$H_{3}=0\oplus0\oplus0\oplus K\oplus K\oplus K,$ \\
$H_{4}=0\oplus0\oplus0\oplus K\oplus\{(y,Sy)\in K^2 \ | \ y\in K\},$
$H_{5}=0\oplus0\oplus0\oplus K\oplus0\oplus0,$ \\
$H_{1^{^{\prime}}}=K\oplus K\oplus\{(x,y,x,y)\in K^4 \ | \ x,y\in K\},$
$H_{2^{^{\prime}}}=K\oplus K\oplus0\oplus0\oplus0\oplus0,$ \\
$H_{1^{^{\prime\prime}}}=\{(y,z,x,0,y,z)\in K^6 \ | \ x,y,z\in K\}.$\\
For any arrow $\alpha \in E$, let  
$f_{\alpha} : H_{s(\alpha)} \rightarrow H_{r(\alpha)}$ be 
the canonical inclusion map. We shall show that 
the Hilbert representation $(H,f)$ is indecomposable. 
Take $T = (T_v)_{v \in V} \in Idem (H,f)$. 
Since $T \in End (H,f)$ and any arrow is represented 
by the inclusion map, we have $T_0 x = T_vx$ for any 
$v \in V $ and 
any  $x \in H_v$. In particular, 
$T_0 H_v \subset H_v$.  
Since $T_{0}$ preserves 
subspaces 
$H_{2^{^{\prime}}}=K\oplus K\oplus0\oplus0\oplus0\oplus0$, 
$H_{2}=0\oplus0\oplus K\oplus K\oplus K\oplus K$, 
$T_0$ has a form such that 
\[
T_{0}=\left(
\begin{array}
[c]{cccccc}%
\ast & \ast & 0 & 0 & 0 & 0\\
\ast & \ast & 0 & 0 & 0 & 0\\
0 & 0 & \ast & \ast & \ast & \ast\\
0 & 0 & \ast & \ast & \ast & \ast\\
0 & 0 & \ast & \ast & \ast & \ast\\
0 & 0 & \ast & \ast & \ast & \ast
\end{array}
\right) 
=
\left(
\begin{array}
[c]{cc}%
A & 0\\
0 & B
\end{array}
\right) , 
\]
for some $A \in B(K \oplus K)$ and $B \in 
B(K \oplus K \oplus K\oplus K)$. 

Moreover 
$H_{1^{^{\prime\prime}}}\cap H_{2}=0\oplus0\oplus K\oplus
0\oplus0\oplus0$ and 
$H_{3}=0\oplus0\oplus0\oplus K\oplus K\oplus K,$
are invariant under $T_0$. Furthermore 
$H_{5}=0\oplus0\oplus0\oplus K\oplus0\oplus0$
and $T_0(H_5) \subset H_5$.
Therefore 
$T_{0}$ is written as 
\[
T_{0}
=
\left(
\begin{array}
[c]{cc}%
A & 0\\
0 & B
\end{array}
\right)
=\left(
\begin{array}
[c]{cccccc}%
a & b & 0 & 0 & 0 & 0\\
c & d & 0 & 0 & 0 & 0\\
0 & 0 & e & 0 & 0 & 0\\
0 & 0 & 0 & f & g & h\\
0 & 0 & 0 & 0 & i & j\\
0 & 0 & 0 & 0 & k & l
\end{array}
\right) , 
\]
for some 
$a,b,c,d,e,f,g,h,i,j,k,l \in B(K)$. 

Since  $H_{1^{^{\prime}}}\cap H_{3}=0\oplus0\oplus0\oplus
\{(y,0,y) \in K^4 \ | \ y\in K\}$ 
is invariant under $T_0$, 
for any $y \in K$, there exists $y' \in K$ such that 
\[
B\left(
\begin{array}
[c]{c}%
0\\
y\\
0\\
y
\end{array}
\right)  =\left(
\begin{array}
[c]{cccc}%
e & 0 & 0 & 0\\
0 & f & g & h\\
0 & 0 & i & j\\
0 & 0 & k & l
\end{array}
\right)  \left(
\begin{array}
[c]{c}%
0\\
y\\
0\\
y
\end{array}
\right)  =\left(
\begin{array}
[c]{c}%
0\\
fy+hy\\
jy\\
ly
\end{array}
\right)  =\left(
\begin{array}
[c]{c}%
0\\
y^{^{\prime}}\\
0\\
y^{^{\prime}}%
\end{array}
\right)  .
\]
Therefore $ f+h  = l$ and $j = 0$. 
Next consider 
$H_{1^{^{\prime}}}\cap H_{2}=0\oplus0\oplus\{(x,y,x,y);x,y\in K\}.$
Since  $\ H_{1^{^{\prime}}}\cap H_{2}$ is invariant under $T_0$, 
for any $x,y \in K$ there exist $x', y' \in K$ such that 
\[
B \left(
\begin{array}
[c]{c}%
x\\
y\\
x\\
y
\end{array}
\right)  =\left(
\begin{array}
[c]{cccc}%
e & 0 & 0 & 0\\
0 & f & g & h\\
0 & 0 & i & 0\\
0 & 0 & k & l
\end{array}
\right)  \left(
\begin{array}
[c]{c}%
x\\
y\\
x\\
y
\end{array}
\right)  =\left(
\begin{array}
[c]{c}%
ex\\
fy+gx+hy\\
ix\\
kx+ly
\end{array}
\right)  =\left(
\begin{array}
[c]{c}%
x'\\
y'\\
x'\\
y'%
\end{array}
\right)
\]
Putting $y = 0$, we have 
\[
ex=x' = ix, \ gx=y' = kx \ \ \text{ for any } x \in K. 
\]
Hence $e=i$ and $g=k$. 

Letting $x = 0$, we have $fy+hy=y' = ly$ for any $y \in K$. 
Hence $f+h=l.$

Since $T_0$ preserves 
$H_{2^{^{\prime}}}\cap H_{1}
=\{(x,x)\in K^2 \ | \ x\in K\}\oplus0\oplus0\oplus0\oplus0, $
for any $x \in K, $there exists $x' \in K$ such that
\[
A\left(
\begin{array}
[c]{c}%
x\\
x
\end{array}
\right)  =\left(
\begin{array}
[c]{cc}%
a & b\\
c & d
\end{array}
\right)  \left(
\begin{array}
[c]{c}%
x\\
x
\end{array}
\right)  =\left(
\begin{array}
[c]{c}%
ax+bx\\
cx+dx
\end{array}
\right)  =\left(
\begin{array}
[c]{c}%
x'\\
x'%
\end{array}
\right)  .
\]
Hence $ax+bx=cx+dx,$ for any $x \in K$, so that $a+b=c+d.$

Furthermore 
$\ H_{1^{^{\prime\prime}}}
=\{(y,z,x,0,y,z)\in K^6 \ | \ x,y,z\in K\}$ 
is invariant under  $T_0$. Therefore for any 
$x,y,z \in K$ there exist $x',y',z' \in K$ satisfying 
\[
\left(
\begin{array}
[c]{cccccc}%
a & b & 0 & 0 & 0 & 0\\
c & d & 0 & 0 & 0 & 0\\
0 & 0 & e & 0 & 0 & 0\\
0 & 0 & 0 & f & g & h\\
0 & 0 & 0 & 0 & e & 0\\
0 & 0 & 0 & 0 & g & l
\end{array}
\right)  \left(
\begin{array}
[c]{c}%
y\\
z\\
x\\
0\\
y\\
z
\end{array}
\right)  
=\left(
\begin{array}
[c]{c}%
ay+bz\\
cy+dz\\
ex\\
gy+hz\\
ey\\
gy+lz
\end{array}
\right)  
=\left(
\begin{array}
[c]{c}%
y'\\
z'\\
x'\\
0\\
y'\\
z'%
\end{array}
\right)  .
\]
Put $x=z=0$. Then for any $y\in K$, we have 
$ay = y' = ey$, $cy = z' = gy$ and $gy = 0$. 
Hence we have $a=e$ and $c = g = 0$.  

Letting $x = y = 0$, for any $z \in K$ we have 
$bz = y' = 0$, $dz = z' = lz$ and $hz = 0$. 
Therefore $b = 0$, $d = l$ and $h = 0$. 
Combining these with $f+h=l$ and $a+b=c+d$, we have
$a = d$ and $f = l = d$. 
Thus $T_0$ is a block diagonal such that 
\[
T_0 = a \oplus a \oplus a \oplus a \oplus a 
      \oplus a \oplus a \oplus a .
\]
Since $T_{0}$ is idempotent, $a$ is also idempotent.

Finally consider that $\ H_{4}
=0\oplus0\oplus0\oplus K\oplus\{(y,Sy)\in K^2 \ | \ y\in K\}$ 
is invariant under $T_0$. 
Then for any $x,y \in K$, there exist $x',y' \in K$ such that 
\[
T_{0}(0,0,0,x,y,Sy)
=(0,0,0,ax,ay,aSy)=(0,0,0,x',y',Sy').
\]
Hence $aSy = Sy' = Say$, so that $aS=Sa$. 
Since $S$ is a unilateral shift and $a$ is idempotent,
we have $a=0$ or $a=I.$
This implies that $T_{0}=0$ or $T_{0}=I.$
Since for any $v \in V$ and any  $x \in H_v$ 
$T_0 x = T_vx$, 
we have $T_v=0$ or $T_v=I$ simultaneously. 
Thus $T = (T_v)_{v \in V} = 0$ or $T = I$, that is, 
$Idem(H,f)=\{0,I\}.$
Therefore $(H,f)$ is indecomposable.
\end{proof}

\smallskip
\noindent
{Remark.} In many cases of our construction of indecomposable, 
infinite-dimensional representations, we can replace a 
unilateral shift $S$ by any strongly irreducible operator.  

We shall show that the existence of indecomposable, 
infinite-dimensional representations does not depend on the 
choice of the orientation of quivers. 
Suppose that two finite, connected quivers 
$\Gamma$  and $\Gamma '$ have the same 
underlying undirected graph and one of them, 
say $\Gamma$, has an infinite-dimensional, 
indecomposable, Hilbert representation. 
We need to prove that another quiver 
$\Gamma '$ also has 
 an infinite-dimensional, 
indecomposable, Hilbert representation. 
Reflection functors are useful to show it. 
But we need to check the co-closedness at 
a source.  We introduce a certain nice class of 
Hilbert representations such that co-closedness is 
easily checked and preserved under reflection functors 
at any source. 

\bigskip
\noindent
{\bf Definition}
Let $\Gamma$ be a quiver whose underlying undirected graph 
is Dynkin diagram $A_n$. We count the arrows from the left as 
$\alpha_k : s(\alpha_k) \rightarrow r(\alpha_k), 
\ (k = 1, \dots, n-1)$.    
Let $(H,f)$ be a Hilbert representation of $\Gamma$. We denote 
$f_{\alpha_k}$ by $f_k$ for short. For example, 
\[
\circ_{H_1} \overset{f_1}{\longleftarrow}
 \circ_{H_2} \overset{f_2}\longrightarrow \circ_{H_3} 
\overset{f_3}{\longleftarrow} 
\circ_{H_4} \overset{f_4}\longrightarrow \circ_{H_5} 
\overset{f_5}\longrightarrow \circ_{H_6}
\]
We say that $(H,f)$ is {\it positive-unitary diagonal} 
if there exist $m \in \mathbb N$ and 
orthogonal decompositions (admitting zero components) 
of Hilbert spaces 
\[
 H_k = \oplus _{i = 1}^m H_{k,i} \ \ \ (k = 1, \dots, n)
\]
and decompositions of operators 
\[
f_k =  \oplus _{i = 1}^m f_{k,i}:  
\oplus _{i = 1}^m H_{s(\alpha_k),i} 
\rightarrow \oplus _{i = 1}^m H_{r(\alpha_k),i} \  
\ \ \ (k = 1, \dots, n), 
\]
such that each 
$f_{k,i} : H_{s(\alpha_k),i} \rightarrow H_{r(\alpha_k),i}$ 
is written as  $f_{k,i} = 0$ or 
$f_{k,i} = \lambda _{k,i} u_{k,i}$ for 
some positive scalar $\lambda _{k,i}$ and onto unitary 
$u_{k,i} \in B(H_{s(\alpha_k),i},H_{r(\alpha_k),i})$. 

It is easy to see that if $(H,f)$ is positive-unitary diagonal, 
then $\Phi^*(H,f)$ is also positive-unitary diagonal.

\bigskip
\noindent
{\bf Example.} (Inclusions of subspaces)
Consider the following quiver $\Gamma$ : 
\[
\circ_1 \overset{\alpha_1}{\longrightarrow}
 \circ_2 \overset{\alpha_2}\longrightarrow \circ_3
\]

Let $H_3$ be a Hilbert space and 
$H_1 \subset H_2  \subset H_3$
inclusions of subspaces. Define a Hilbert 
representation $(H,f)$ of $\Gamma$ by $H = (H_i)_{i= 1,2,3}$ and 
canonical inclusion maps $f_i = f_{\alpha _i} : H_i \rightarrow H_{i+1}$ 
for $i = 1,2$. Then $(H,f)$ is positive-unitary diagonal. 
In fact, define 
\[
K_1 = H_1,\ K_2 = H_2 \cap H_1^{\perp}, \ 
K_3 = H_3 \cap H_2^{\perp}.
\]
Consider orthogonal decompositions 
$H_k = \oplus _{i = 1}^3 H_{k,i} \ (k = 1,2,3)$ 
by 
\[
H_1 = K_1 \oplus 0 \oplus 0, 
H_2 = K_1 \oplus K_2 \oplus 0 \text{ and } 
H_3 = K_1 \oplus K_2 \oplus K_3. 
\]  
Then $f_1= I \oplus  0 \oplus 0$ and 
$f_2 = I \oplus  I \oplus 0$. Hence $(H,f)$ 
is positive-unitary diagonal. It is trivial that 
the example can be extended to the case of 
inclusion of $n$ subspaces. 

\begin{lemma}
Let $\Gamma$ be a quiver whose underlying undirected graph 
is Dynkin diagram $A_n$ and 
$(H,f)$ be a Hilbert representation of $\Gamma$. 
Assume that $(H,f)$ is positive-unitary diagonal. Then 
$(H,f)$ is closed at any sink of $\Gamma$ and 
co-closed at any source of $\Gamma$. 
\label{lemma:positive-unitary-coslosed} 

\end{lemma}
\begin{proof}Let $v$ be a sink of $\Gamma$. 
Then 
$\sum_{\alpha \in E^v} \Im f_{\alpha}$
is a finite sum of some of orthogonal subspaces 
$\{H_{v,i} \ | i \}$ of $H_v$ which correspond to 
the images of positive times unitaries in the direct 
component of $f_{\alpha}$. 
Hence it is a closed subspace of $H_v$. Therefore 
$(H,f)$ is closed at $v$.  Similarly $(H,f)$
co-closed at any source of $\Gamma$.
\end{proof}

\begin{prop}
Let $\Gamma$ be a quiver whose underlying undirected graph 
is Dynkin diagram $A_n$ and 
$(H,f)$ be a Hilbert representation of $\Gamma$. 
Let $v$ be a source of $\Gamma$. 
Assume that $(H,f)$ is positive-unitary diagonal. 
Then $\Phi_v^{-}(H,f)$ is also positive-unitary diagonal. 
\label{prop:positive-unitary-iteration} 

\end{prop}
\begin{proof}
If $(H,f) \cong (H',f') \oplus (H'',f'')$, 
then $\Phi_v^{-}(H,f) 
\cong \Phi_v^{-}(H',f') \oplus \Phi_v^{-}(H'',f'')$. 
Therefore $H_k^- = \oplus _{i=1}^m H_{k,i}^- $. 
Hence it is enough to consider orthogonal components. 
We may and do examine locally the following cases:

\noindent
(Case 1): A Hilbert representation $(H,f)$ is given by 
\[
\circ_{H_1} \overset{T_1}{\longleftarrow}
 \circ_{H_0} \overset{T_2}\longrightarrow \circ_{H_2} 
\]
with $T_1 = \lambda_1U_1$ and $T_2 = \lambda_2U_2$ 
for some positive scalars $\lambda_1,  \lambda_2$ and 
onto unitaries $U_1, U_2$. 
Put $(H^-,f^-) = \Phi_0^{-}(H,f)$: 
\[
\circ_{H_1} \overset{T_1^-}{\longrightarrow}
 \circ_{H_0^-} \overset{T_2^-}\longleftarrow \circ_{H_2} 
\]
Then $(a,b)\in  H_1 \oplus H_2$ is in 
$H_0^- = (\Im \hat{h}_0)^{\perp}$ if and only if 
$((a,b) \ | \  (T_1z,T_2z)) = 0$ for any $z \in H_0$, so that 
$T_1^*a + T_2^*b = 0$. Hence 
\begin{align*}
& H_0^- = \{(a,-\lambda_1\lambda_2^{-1}U_2U_1^*a) 
\in H_1 \oplus H_2  \ | \ a \in H_1 \} \\
& = \{(-\lambda_1^{-1}\lambda_2U_1U_2^*b,b) \in H_1 \oplus H_2
\ | \ b \in H_2 \}. 
\end{align*}
Solving 
\[
(x,0) = 
(a,-\lambda_1\lambda_2^{-1}U_2U_1^*a) + (\lambda_1U_1z,\lambda_2U_2z) 
\in H_0^- \oplus \Im \hat{h}_0 , 
\]
we have 
\[
T_1^- x = (\frac{\lambda_2^2}{\lambda_1^2 + \lambda_2^2}x, 
 -\frac{\lambda_1\lambda_2}{\lambda_1^2 + \lambda_2^2}
 U_2U_1^*x)  \text{ for } x \in H_1 . 
\]
Similarly we have 
\[
T_2^- y 
= (-\frac{\lambda_1\lambda_2}{\lambda_1^2 + \lambda_2^2}U_1U_2^*y, 
\frac{\lambda_2^2}{\lambda_1^2 + \lambda_2^2}y) 
\text{ for } y \in H_2  . 
\]
Let 
$\lambda_1^- := \sqrt{(\frac{\lambda_2^2}{\lambda_1^2 + \lambda_2^2})^2 
+ (\frac{\lambda_1\lambda_2}{\lambda_1^2 + \lambda_2^2} )^2  } > 0$ 
and $U_1^-:= (\lambda_1^-)^{-1}T_1^-$. Then 
$U_1^-$ is an onto unitary and $T_1^- = \lambda_1^-U_1^-$. 
Similarly $T_2^-$ is  a positive scalar times unitary. 

\noindent
(Case 2):A Hilbert representation $(H,f)$ is given by 
\[
\circ_{H_1} \overset{T_1}{\longleftarrow}
 \circ_{H_0} \overset{T_2}\longrightarrow \circ_{H_2} 
\]
with $T_1 = 0$ and $T_2 = 0$

Then it is easy to see that $H_0^- = H_1 \oplus H_2$, 
$T_1^-$and $T_2^-$ are canonical inclusions: 
$T_1^-x = (x,0)\in H_1 \oplus H_2$ for $x \in H_1$ 
and $T_2^-y = (0,y) \in H_1 \oplus H_2$ for $y \in H_2$. 
We may write that 
$T_1^- = I \oplus 0: H_1 \oplus 0 \rightarrow H_1 \oplus H_2$ 
and 
$T_2^- = 0 \oplus I: 0 \oplus H_2 \rightarrow H_1 \oplus H_2$ .
Hence $(H^-,f^-)$ is positive-unitary diagonal. 
 
\noindent
(Case 3):A Hilbert representation $(H,f)$ is given by 
\[
\circ_{H_1} \overset{T_1}{\longleftarrow}
 \circ_{H_0} \overset{T_2}\longrightarrow \circ_{H_2} 
\]
with $T_1 = \lambda_1U_1$ and $T_2 = 0$ 
for some positive scalar $\lambda_1$ and 
onto unitary $U_1$. 

Then we see that $H_0^- = 0 \oplus H_2$, 
$T_1^- = 0$ and 
$T_2^-y = (0,y) \in 0\oplus H_2$ for $y \in H_2$.
Hence $(H^-,f^-)$ is positive-unitary diagonal. 

\noindent
(Case 4):A Hilbert representation $(H,f)$ is given by 
\[
\circ_{H_0} \overset{T_1}{\longrightarrow}
 \circ_{H_1}
\]
with $T_1 = \lambda_1U_1$ 
for some positive scalar $\lambda_1$ and 
onto unitary $U_1$. 
Put $(H^-,f^-) = \Phi_0^{-}(H,f)$: 
\[
\circ_{H_0^-} \overset{T_1^-}{\longleftarrow}
 \circ_{H_1}  
\]
Then we see that $H_0^- = 0$ and 
$T_1^- = 0$. 

\noindent
(Case 5):A Hilbert representation $(H,f)$ is given by 
\[
\circ_{H_0} \overset{T_1}{\longrightarrow}
 \circ_{H_1}
\]
with $T_1 = 0$.  

Then we have  that
$H_0^- = H_1$ and 
$T_1^- = I : H_1 \rightarrow H_1 = H_0^-$. 
\end{proof}

We shall show that we can change the orientation of Dynkin 
diagram $A_n$ using only  the iteration of 
$\sigma_v^-$ at sources $v$ except the right end.

\begin{lemma}Let $\Gamma _0$  and $\Gamma$ be quivers
 whose underlying undirected graphs
are the same Dynkin diagram $A_n$ for $n \geq 2$. We assume that 
$\Gamma _0$ is the following: 
\[
\circ_1 \longrightarrow \circ_2 \longrightarrow \circ_3 \dots 
\circ_{n-1} \longrightarrow \circ_n
 \]
Then there exists a sequence $v_1, \dots, v_m$ of vertices in 
$\Gamma _0$ such that 
\begin{itemize}
\item [(1)] for each $k= 1, \dots, m$, $v_k$ is a source in 
$\sigma^-_{v_{k-1}} \dots \sigma^-_{v_2} \sigma^-_{v_1}(\Gamma _0)$,
\item [(2)] 
$\sigma^-_{v_m} \dots \sigma^-_{v_2} \sigma^-_{v_1}(\Gamma _0) = 
\Gamma$, 
\item[(3)]for each $k= 1, \dots, m$, $v_k \not= n$. 
\end{itemize}
\label{lemma:orientation-change} 
\end{lemma}
\begin{proof}
The proof is by induction on the number $n$ of vertices. 
Let $n = 2$. Since 
$\sigma^-_1(\circ_1 \longrightarrow \circ_2) = \circ_1 \longleftarrow \circ_2$, the statement holds. Assume that the statement holds for $n-1$. If $\Gamma$ 
has an arrow $\circ_{n-1} \longrightarrow \circ_n$, then we can directly 
apply the assumption of the induction. If $\Gamma$ 
has an arrow $\circ_{n-2} \longrightarrow \circ_{n-1} \longleftarrow \circ_n$, 
replace only this part by 
$\circ_{n-2} \longleftarrow \circ_{n-1} \longrightarrow \circ_n$ to get 
$\Gamma'$.  Then $n-1$ is a source of $\Gamma'$, and 
$\sigma^-_{n-1}(\Gamma') = \Gamma$. Applying the induction assumption 
for  $\Gamma'$, we can construct the desired iteration. 
Consider the case that $\Gamma$ 
has an arrow $\circ_{n-2} \longleftarrow \circ_{n-1} \longleftarrow \circ_n$.
If there exist a vertex $u$ such that $\circ_{u-1} \longrightarrow \circ_{u}$ 
and $\circ_{k} \longleftarrow \circ_{k+1}$ for $k = u, \dots , n-1$, then 
define a new quiver $\Gamma''$ by putting 
$\circ_{u-1}\longleftarrow \circ_{u}$, 
$\circ_{n-1}\longrightarrow \circ_{n}$ and 
other arrows unchanged with $\Gamma$. By the induction assumption, 
there exists a sequence $v_1, \dots, v_m$ of vertices in 
$\Gamma _0$ such that 
$\sigma^-_{v_m} \dots \sigma^-_{v_2} \sigma^-_{v_1}(\Gamma _0) = 
\Gamma ''$ and, for each $k= 1, \dots, m$, $v_k \not= n$ and  $v_k \not= n-1$. 
Then 
\[
\sigma^-_{u} \sigma^-_{u+1} \dots \sigma^-_{n-2} \sigma^-_{n-1}
\sigma^-_{v_m} \dots \sigma^-_{v_2} \sigma^-_{v_1}(\Gamma _0) 
= \Gamma . 
\]

If all the arrows between $1$ and $n$ are of the form 
$\circ_{k} \longleftarrow \circ_{k+1}$ for $k = 1, \dots , n-1$, 
then 
$\sigma^-_{n-1} \dots \sigma^-_{2} \sigma^-_{1}(\Gamma _0) = \Gamma$.
\end{proof}

\begin{lemma}
Let $\Gamma = (V,E,s,r)$  and 
$\Gamma ' = (V',E',s',r')$ be finite, connected quivers 
and $\Gamma '$ contains $\Gamma$ as a subgraph, that is, 
$V \subset V'$, $E \subset E'$, 
$s = s'|_E$ and $r = r'|_E$. 
If there exists an infinite-dimensional, indecomposable, 
 Hilbert representation of $\Gamma$, then 
there exists an infinite-dimensional, indecomposable, 
 Hilbert representation of $\Gamma '$. 
\label{lemma:subgraph} 
\end{lemma}
\begin{proof}
Let $(H,f)$ be an infinite-dimensional, indecomposable, 
 Hilbert representation of $\Gamma$. Define $H'_v = H_v$ 
for $v \in V$ and $H'_v = 0$ for $v \in V' \setminus V$.  
We put $f'_{\alpha} = f_{\alpha}$
for $\alpha  \in E$ and $f'_{\alpha} = 0$ for 
$\alpha \in E' \setminus E$. Then it is clear that 
$(H',f')$ is an infinite-dimensional, indecomposable, 
 Hilbert representation of $\Gamma '$.
\end{proof}

We are ready to prove our main theorem. 
\begin{thm}
Let $\Gamma$ be a finite, connected quiver. 
If the underlying undirected graph $|\Gamma|$ contains 
one of the extended Dynkin diagrams  
$\tilde{A_n} \ (n \geq 0), \tilde{D_n} \  (n \geq 4), 
\tilde{E_6},\tilde{E_7}$ and $\tilde{E_8}$, then 
there exists an infinite-dimensional, indecomposable, 
 Hilbert representation of $\Gamma$. 
\end{thm}
\begin{proof}
By Lemma \ref{lemma:subgraph}, we may assume that 
the underlying undirected graph $|\Gamma|$ is exactly 
one of the extended Dynkin diagrams  
$\tilde{A_n} \ (n \geq 0), \tilde{D_n} \  (n \geq 4), 
\tilde{E_6},\tilde{E_7}$ and $\tilde{E_8}$. 

The case of extended Dynkin diagrams  
$\tilde{A_n} \ (n \geq 0)$ was already verified in 
Example 1 and 2 in section 3.  

Next suppose that $|\Gamma|$ is $\tilde{E_6}$.  Let $\Gamma _0$ 
be the quiver of Example 4 in section 3 and  we denote here by
$(H^{(0)},f^{(0)})$ 
the Hilbert representation constructed there. Then  
$|\Gamma _0| =  |\Gamma| = \tilde{E_6}$, but their orientations 
are different in general. 
Three "wings" of $|\Gamma _0|$ 
$2-1-0,\  2'-1'-0, \ 2''-1''-0$ 
are of  $A_3$. 
Applying Lemma  \ref{lemma:orientation-change}  for 
these wings locally, we can 
find a sequence $v_1, \dots, v_m$ of vertices in 
$\Gamma _0$ such that 
\begin{itemize}
\item [(1)] for each $k= 1, \dots, m$, $v_k$ is a source in 
$\sigma^-_{v_{k-1}} \dots \sigma^-_{v_2} \sigma^-_{v_1}(\Gamma _0)$, 
\item [(2)] 
$\sigma^-_{v_m} \dots \sigma^-_{v_2} \sigma^-_{v_1}(\Gamma _0) = \Gamma$, 
\item[(3)]for each $k= 1, \dots, m$, $v_k \not= 0$. 
\end{itemize}
We note that co-closedness of Hilbert representations at 
a source can be checked locally around the source. 
Since the restriction of the representation $(H^{(0)},f^{(0)})$ 
to each "wing" is positive-unitary diagonal and the iteration of 
reflection functors does not move the vertex $0$, we can apply 
Lemma \ref{lemma:positive-unitary-coslosed} and 
Proposition \ref{prop:positive-unitary-iteration} locally 
that $\Phi^-_{v_{k-1}} \dots 
\Phi^-_{v_2} \Phi^-_{v_1}(H^{(0)},f^{(0)})$ is 
co-closed at $v_k$ for $k= 1, \dots, m$. 
Therefore 
Theorem  \ref{thm:preserving-indecomposability} implies that 
$(H,f) :=  \Phi^-_{v_m} \dots 
\Phi^-_{v_2} \Phi^-_{v_1}(H^{(0)},f^{(0)})$
is the desired indecomposable, 
 Hilbert representation of $\Gamma$. 
Since the particular Hilbert space $H_0^{(0)}$ associated with 
the vertex $0$ is infinite dimensional and remains unchanged 
under the iteration of the reflection functors above, 
$(H,f)$ is infinite dimensional. 

The case that the $|\Gamma|$ is $\tilde{E_7}$ or $\tilde{E_8}$ is 
shown similarly if we apply iteration of reflection functors 
on the representations in Lemma \ref{lemma:E7} 
or Lemma \ref{lemma:E8}. 

Finally consider the case that the $|\Gamma|$ is $\tilde{D_n}$.  
Let $\Gamma _0$ 
be the quiver of Lemma \ref{lemma:Dn}  and $(H^{(0)},f^{(0)})$ 
the Hilbert representation constructed there. Then  
$|\Gamma _0| =  |\Gamma| = \tilde{D_n}$, but their orientations 
are different in general. Let $\Gamma _1$ be a  
quiver such that $|\Gamma _1| = \tilde{D_n}$ and the orientation 
is as same as $\Gamma $ on the path between 5 and n+1 and 
as same as $\Gamma_0$ on the rest four "wings". Define a Hilbert 
representation $(H^{(1)},f^{(1)})$ of $\Gamma _1$ 
similarly as $(H^{(0)},f^{(0)})$. 
For any arrow $\beta$ in the path between 5 and n+1, 
$f^{(1)}_{\beta} = I$. Hence the same proof as for 
$(H^{(0)},f^{(0)})$ shows that $(H^{(1)},f^{(1)})$ is indecomposable. 
By a certain iteration of reflection functors at a source 
1,2,3 or 4 on $(H^{(1)},f^{(1)})$ yields an 
infinite-dimensional, indecomposable, 
 Hilbert representation of $\Gamma$. Here the co-closedness 
at a source 
1,2,3 or 4 on $(H^{(1)},f^{(1)})$ is easily 
checked, because the map is the canonical inclusion. Thus 
we can apply Theorem  \ref{thm:preserving-indecomposability} 
in this case too. 

\end{proof}

\begin{cor}
Let $\Gamma$ be a finite, connected  quiver. 
If there exists no infinite-dimensional, indecomposable, 
 Hilbert representation of $\Gamma$, then 
the underlying undirected graph $|\Gamma|$ is one of the 
Dynkin diagrams $A_n \ (n \geq 1), D_n \  (n \geq 4), 
E_6, E_7$ and $E_8$. 
\end{cor}
\begin{proof}
It directly follows from a well known fact that 
if the underlying undirected graph $|\Gamma|$ contains no 
extended Dynkin diagrams, then $|\Gamma|$ is one of the 
 Dynkin diagrams. 
\end{proof}

\noindent
{\bf Remark.}
We have not yet  proved the converse. 
In fact if the converse were true, then a long standing 
problem on  transitive lattices of subspaces of 
Hilbert spaces would be settled.  Recall that 
Halmos initiated the study of transitive lattices and gave an example of 
transitive lattice consisting of seven subspaces in \cite{Ha}. 
Harrison-Radjavi-Rosenthal \cite{HRR} constructed a transitive lattice  
consisting of six subspaces using the graph of an unbounded closed operator. 
Hadwin-Longstaff-Rosenthal found a transitive lattice of five  non-closed 
linear subspaces in \cite{HLR}.     
Any finite transitive lattice which consists of  $n$ subspaces of a Hilbert 
space $H$  
gives an indecomposable system of $n-2$ subspaces by withdrawing 
$0$ and $H$.  It is still unknown whether or not 
there exists a transitive lattice 
consisting of five subspaces.  Therefore it is also an interesting  
problem to know whether there exists an indecomposable system of 
three subspaces in an infinite-dimensional Hilbert space.  
The problem can be rephrased as whether there exists an 
indecomposable representation of a certain quiver whose underlying 
undirected graph is $D_4$
in an infinite-dimensional Hilbert space. 

We have a partial evidence for a certain quiver whose underlying 
undirected graph is $A_n$.  We prepare an elementary lemma.  Let 
$H$ be a Hilbert space. For $a,b \in H$ we denote by $\theta _{a,b}$ 
a rank one operator on $H$ such that $\theta _{a,b}(x) = (x | b)a$ 
for $x \in H$. Then $\theta _{a,b}^2 = \theta _{a,b}$ if and only if 
$(a|b) = 1$ or $a = 0$ or $b = 0$. Moreover if $\dim H \geq 2$ and 
$(a|b) = 1$, then  $\theta _{a,b}$ is an idempotent such that 
$\theta _{a,b} \not= 0$ and $\theta _{a,b} \not= I$. 

\begin{lemma}
Let $H_1$ and $H_2$ be Hilbert spaces and 
$T: H_1 \rightarrow H_2$ a bounded operator. Take 
$a,b \in H_1$ and $c,d \in H_2$. Suppose that 
there exists a scalar $\lambda $ such that $Ta = \lambda c$ and 
$T^*d = \overline{\lambda} b$. Then 
$T\theta _{a,b} = \theta _{c,d} T$
\label{lemma:rank-one}
\end{lemma}
\begin{proof}
\[
T\theta _{a,b} = \theta _{Ta,b} = \theta _{\lambda c,b} 
= \theta _{c,\overline{\lambda} b} = \theta _{c,T^*d} = \theta _{c,d} T. 
\]
\end{proof}

\begin{prop}
Let $\Gamma$ be the following quiver 
 whose underlying undirected graph is $A_n$ for $n \geq 1$:
\[
\circ_1 \overset{\alpha_1}{\longrightarrow} \circ_2 
\overset{\alpha_2}{\longrightarrow} \circ_3 \dots 
\circ_{n-1} 
\overset{\alpha_{n-1}}{\longrightarrow} \circ_n
\]
Then there exists no 
infinite-dimensional, indecomposable, 
 Hilbert representation of $\Gamma$. 
\end{prop}
\begin{proof} The case $n = 1$ is clear by a nontrivial decomposition 
$H_1 = L_1 \oplus K_1$.  
We may assume 
that $n \geq 2$. 
Suppose that there were an infinite-dimensional, indecomposable, 
 Hilbert representation $(H,f)$ of $\Gamma$. Put 
$T_k = f_{\alpha _k} : H_k \rightarrow H_{k+1}$ for $k = 1,\dots , n-1$. 

\noindent
(case 1) Suppose that $T_{n-1}T_{n-2} \dots T_1 \not= 0$. Then 
there exists $a_1 \in H_1$ such that 
$T_{n-1}T_{n-2} \dots T_1 a_1 \not= 0 $ . Consider non-zero vectors  
$a_k = T_{k-1}T_{k-2} \dots T_1 a_1 \in H_k$ for $k = 1,\dots, n$. 
Put $b_n = \|a_n\|^{-2}a_n \in H_n$. Define 
$b_i = T_{i}^*T_{i+1}^* \dots T_{n-1}^*b_n \in H_i$ for $i = 1,2, \dots, n-1$. 
Then 
\[
(a_i |b_i) = (a_i|T_{i}^*T_{i+1}^* \dots T_{n-1}^*b_n) 
= (T_{n-1}T_{n-2} \dots T_i a_i|b_n) 
=  (a_n|b_n) = 1. 
\]
Since $T_ka_k = a_{k+1}$ and $T_k^*b_{k+1} = b_k$, 
the above Lemma \ref{lemma:rank-one} implies that 
$T_k\theta _{a_k,b_k} = \theta _{a_{k+1},b_{k+1}} T_k$ 
for  $k = 1, \dots , n-1$. 
Define the non-zero idempotents $P_k = \theta _{a_k,b_k}$. 
Since $(H,f)$ is infinite dimensional, there exists some 
vertex $m$ such that $H_m$ is infinite dimensional. 
Then $P_m \not= I$. Define $P = (P_k)_k$, 
then $ P \in Idem(H,f)$ and $P \not= O$ and $P \not= I$. 
This contradicts the assumption that $(H,f)$ is indecomposable. 

\noindent
(case 2) Suppose that there exists $r$ such that 
$T_{r-1}T_{r-2} \dots T_1 \not= 0$ and $T_{r}T_{r-1} \dots T_1 = 0$ 
for some $r = 1,\dots, n-1$ and $\dim H_m \geq 2$  for some 
$m = 1, \dots ,r$. 
Then there exists $a_1 \in H_1$ such that 
$T_{r-1}T_{r-2} \dots T_1 a_1 \not= 0 $ . 
Consider non-zero vectors  
$a_k = T_{k-1}T_{k-2} \dots T_1 a_1 \in H_k$ for $k = 1,\dots, r$. 
Put $b_r = \|a_r\|^{-2}a_r \in H_r$. Define 
$b_i = T_{i}^*T_{i+1}^* \dots T_{r-1}^*b_r \in H_i$ for 
$i = 1,2, \dots, r-1$. 
Then we have $T_k\theta _{a_k,b_k} = \theta _{a_{k+1},b_{k+1}} T_k$ 
for  $k = 1, \dots , r-1$ as case 1. 
Define non-zero idempotents $P_k = \theta _{a_k,b_k}$ for 
$k = 1, \dots , r$. Put $P_k = 0$ for 
$k = r+1, \dots , n$. Then 
$T_r\theta _{a_r,b_r} = \theta _{T_ra_r,b_r} = \theta _{0,b_r} =0$ and 
$T_kP_k = P_{k+1} T_k = 0$ for 
$k = r, \dots , n-1$. 
Since 
$dim H_m \geq 2$,  the non-zero idempotent $P_m \not= I$.
 Define $P = (P_k)_k$, 
then $ P \in Idem(H,f)$ and $P \not= O$ and $P \not= I$. 
This is a contradiction. 

\noindent
(case 3) Suppose that there exists $r$ such that 
$T_{r-1}T_{r-2} \dots T_1 \not= 0 $ and $T_{r}T_{r-1} \dots T_1 = 0$ 
for some $r = 1,\dots, n$ and  $\dim H_k = 1$ for 
$k = 1, \dots, r$. Therefore $T_r = 0$. We may put 
$P_k = 0$ for $k = 1, \dots, r$. Then for any $a,b \in H_{r+1}$ and 
$P_{r+1} = \theta _{a,b}$, we have $T_kP_k = P_{k+1} T_k = 0$ for 
$k = 1, \dots , r$. Hence we may choose freely $P_k$ for 
$k = r+1 , ..., n$.  Starting form $H_{r+1}$, we can repeat the argument 
from the beginning.  After finite steps, we can reduce to the situation of case 1 or case 2. And finally we obtain a contradiction.

\end{proof}

\end{document}